\documentclass{article}
\usepackage{geometry}
\usepackage{tikz}
\usepackage{graphicx}
\usepackage{enumerate}
\usepackage{amsmath}
\usepackage{amssymb}
\usepackage{amsthm}
\usepackage{hyperref}

\newcommand{\mprog}[0]{$({\leq}m)$-progression}
\newcommand\tuple[1]{\mathbf{#1}} 

\newcommand{\Pol}[0]{\ensuremath{\mathrm{Pol}}}
\newcommand{\Aut}[0]{\ensuremath{\mathrm{Aut}}}
\newcommand{\End}[0]{\ensuremath{\mathrm{End}}}

\renewcommand{\phi}{\varphi}
\newcommand{\ignore}[1]{}

\newcommand{\N}{\mathbb{N}}
\newcommand{\Z}{\mathbb{Z}}
\newcommand{\Q}{\mathbb{Q}}

\newcommand{\Diff}[1]{\ensuremath{\mathrm{Diff}_{#1}}}
\newcommand{\Dist}[1]{\ensuremath{\mathrm{Dist}_{#1}}}

\DeclareMathOperator{\suc}{succ}
\DeclareMathOperator{\Mod}{Mod}

\DeclareMathOperator{\Csp}{CSP}
\DeclareMathOperator{\Th}{Th}

\DeclareMathOperator{\si}{si}

\newtheorem{theorem}{Theorem}
\newtheorem*{thm:class}{Theorem~\ref{thm:class}}
\newtheorem*{thm:main}{Theorem~\ref{thm:main}}
\newtheorem*{thm:succ}{Theorem~\ref{thm:succ}}
\newtheorem{corollary}{Corollary}
\newtheorem{lemma}{Lemma}
\newtheorem{proposition}{Proposition}
\theoremstyle{definition}
\newtheorem{definition}{Definition}
\theoremstyle{remark}
\newtheorem{example}{Example}
\newtheorem{remark}{Remark}

\begin{document}

\title{Discrete Temporal Constraint Satisfaction Problems\footnote{A subset of the results appeared in the Proceedings of ICALP'15 under the title \emph{Constraint Satisfaction Problems over the Integers with Successor}.}}
\author{Manuel Bodirsky\thanks{The first and third authors have received funding from the European Research Council under the European Community's Seventh Framework Programme (FP7/2007-2013 Grant Agreement no. 257039), and the DFG (Project number 622397).}\\Institut f\"ur Algebra, TU Dresden, Dresden
\and Barnaby Martin\thanks{This author was supported by EPSRC grant EP/L005654/1.}\\School of Science and Technology, Middlesex University, London
\and Antoine Mottet$^\dagger$\thanks{This author is supported by DFG Graduiertenkolleg 1763 (QuantLA).}\\Institut f\"ur Algebra, TU Dresden, Dresden
}

\maketitle

%

\begin{abstract}
A \emph{discrete temporal constraint satisfaction problem} is a constraint satisfaction problem (CSP) whose constraint language consists of relations that are first-order definable over $({\mathbb Z};<)$.
Our main result says that every discrete temporal CSP is in Ptime or NP-complete,
unless it can be formulated as a finite domain CSP in which case the computational complexity is not known in general. 
\end{abstract}

\tableofcontents

\section{Introduction}
\begin{flushright}
\emph{``Die ganzen Zahlen hat der liebe Gott gemacht, alles andere ist Menschenwerk.''}\footnote{\emph{``God made the integers, all the rest is the work of man.''} Quoted in \emph{Philosophies of Mathematics}, page 13, by Alexander George, Daniel J. Velleman, Philosophy, 2002.} Leopold Kronecker
\end{flushright}
A \emph{constraint satisfaction problem} is a computational problem where the input consists of a finite set of variables and a finite set of constraints, and where the question is whether there exists a mapping from the variables to 
some fixed domain such that all the constraints are satisfied. 
When the domain is finite, and arbitrary constraints are permitted in the input, the CSP is NP-complete. 
However, when only constraints for a restricted set of relations are allowed in the input, it might be possible to solve the CSP in polynomial time. 
The set of relations that is allowed to formulate the constraints in the input is often called the \emph{constraint language}. The question which constraint
languages give rise to polynomial-time solvable CSPs
has been the topic of intensive research over the past years. It has been conjectured by Feder and Vardi~\cite{FederVardi} that CSPs for constraint languages over finite domains have a complexity dichotomy: they are in Ptime or NP-complete. 

A famous CSP over an infinite domain is \emph{feasibility of linear inequalities over the integers}.
It is of great importance in practice and theory of computing, and NP-complete. 
In order to obtain a systematic understanding
of polynomial-time solvable restrictions and variations of this problem, 
Jonsson and L\"o\"ow~\cite{JonssonLoeoew}  proposed to study the class of
CSPs where the constraint language $\Gamma$ 
is definable in 
\emph{Presburger arithmetic}; that is,
it consists
of relations that have a first-order definition over $({\mathbb Z};\leq,+)$. 
Equivalently, each relation $R(x_1,\dots,x_n)$ in $\Gamma$ can be
defined by a disjunction of conjunctions of the atomic formulas of the form $p \leq 0$ where
 $p$ is a linear polynomial with integer coefficients and  variables from $\{x_1,\dots,x_n\}$. 
The constraint satisfaction problem for $\Gamma$, denoted by $\Csp(\Gamma)$, is the problem of deciding whether a
given conjunction of formulas of the form $R(y_1,\dots,y_n)$, for some $n$-ary $R$ from $\Gamma$, is satisfiable in $\Gamma$. 
By appropriately choosing such a constraint language $\Gamma$, a great
variety of problems over the integers can be formulated as $\Csp(\Gamma)$.
Several constraint languages $\Gamma$ over the integers are known where the CSP can be solved in polynomial time.
However, a complete complexity classification for the CSPs of Jonsson-L\"o\"ow languages appears to be a very ambitious goal. 

One of the most basic classes of constraint languages that falls into the framework of Jonsson and L\"o\"ow
is the class of \emph{distance constraint satisfaction problems}~\cite{BodDalMarMotPin}. 
A distance constraint satisfaction problem is a CSP for a constraint language over the integers whose relations
have a first-order definition over $(\Z;\suc)$ where $\suc$ is the successor function.
It has been shown previously that distance CSPs for constraint languages whose relations have \emph{bounded Gaifman degree} are either NP-complete, or in Ptime,
or can also be formulated with a constraint language over a finite domain~\cite{BodDalMarMotPin}. 

Another class of problems which can be expressed as Jonsson-L\"o\"ow constraint satisfaction problems is the class of \emph{temporal CSPs}~\cite{tcsps-journal}.
This is the class of problems whose constraint languages are structures which are definable over $(\Q;<)$.
While the order of the rationals is not isomorphic to the order of integers because of its density,
this density is not witnessed by \emph{finite structures}: for any finite substructure of $(\Q;<)$, one can find
a substructure of $(\Z;<)$ that is order-isomorphic to it.
It follows that for every first-order reduct $\Gamma$ of $(\Q;<)$, there exists a reduct $\Delta$ of $(\Z;<)$
such that $\Gamma$ and $\Delta$ have the same constraint satisfaction problem.
In the present paper, we prove that the class of \emph{discrete temporal CSPs}, that is, the constraint satisfaction problems
whose constraint language is first-order definable over $(\Z;<)$, exhibits a P/NP-complete dichotomy (modulo the Feder-Vardi conjecture for finite-domain CSPs).
Our result thus properly extends the results mentioned above for locally finite distance CSPs and temporal CSPs. For the case of distance CSPs,
our approach replaces the proof of~\cite{BodDalMarMotPin}.
However, a cornerstone of our proof is the characterization of those problems that are discrete temporal CSPs but that are not temporal CSPs.

Our proof relies on the so-called universal-algebraic approach; this is the first-time that this approach has been used for constraint languages that are not finite or countably infinite \emph{$\omega$-categorical}.
The central insight of the universal-algebraic approach to constraint satisfaction is that the computational complexity of a CSP is captured by the set of \emph{polymorphisms} of the constraint language. One of the ideas of the present paper is that in order to use polymorphisms when the constraint language is not $\omega$-categorical, we have to pass to the countable saturated model of $(\Z;<)$.  
The relevance of saturated models for the universal-algebraic approach has already been pointed out in joint work of the first two authors with Martin Hils~\cite{BodHilsMartin-Journal},
but this is the first time that this perspective has been used to perform complexity classification for a large class of concrete computational problems. 

The formal definitions of CSPs and discrete temporal CSPs can be found in Section~\ref{sect:prelims}. 
The border between discrete temporal CSPs in Ptime and NP-complete discrete temporal CSPs can be most elegantly stated using the terminology of the mentioned universal-algebraic approach to constraint satisfaction. This is why we first 
give a brief introduction to this approach in Section~\ref{sect:algebra}, and only then
give the technical description 
of our main result in Section~\ref{sect:results}. 
Section~\ref{sect:succ} gives
a classification of discrete temporal constraint languages
that might be of independent interest;
this classification is the basis of our classification of the complexity
of discrete temporal CSPs.
Our algorithmic results can be found in Section~\ref{sect:algs}.
Finally, we put all the results together to prove our main result in Section~\ref{sect:classification}. 
We discuss our result and promising future research questions in Section~\ref{sect:discussion}.

\section{Discrete Temporal CSPs}
\label{sect:prelims}
Let $\Gamma$ be a structure
with a finite relational signature $\tau$. 
When $R$ is a relation symbol from $\tau$,
we write $R^\Gamma$ for the relation it denotes in the structure $\Gamma$.

A \emph{$\tau$-formula} is a first-order formula
built from the relations from $\tau$, and equality. A $\tau$-formula is 
\emph{primitive positive (pp)}
if it is of the form $\exists x_1,\dots,x_k (\psi_1 \wedge \cdots \wedge \psi_m)$ where each
$\psi_i$ is an atomic $\tau$-formula. 
\emph{Sentences} are formulas without free variables. 

We define the constraint satisfaction problem for $\Gamma$ as follows.

\begin{definition}[$\Csp(\Gamma)$]
The \emph{constraint satisfaction problem for $\Gamma$} is the following computational problem. \\
{\bf Input:} A primitive positive $\tau$-sentence $\Phi$.\\
{\bf Question:} $\Gamma \models \Phi$?
\end{definition}

The structure $\Gamma$ will also be called the
\emph{constraint language} of $\Csp(\Gamma)$. 
A relational structure $\Gamma$ is a \emph{reduct} of a structure $\Delta$ if it has the same domain 
as $\Delta$ and for every relation $R^\Gamma$ of arity $k$ is \emph{first-order definable} over $\Delta$, that is,
there exists a first-order formula $\phi$ in the signature of $\Delta$ with $k$ free variables such that for
all elements $u_1,\dots,u_k$ of $\Gamma$ we have  $R^\Gamma(u_1,\dots,u_k) \Leftrightarrow \Delta \models \phi(u_1,\dots,u_k)$.

\begin{definition}[Discrete Temporal CSP]
A \emph{discrete temporal CSP} is a constraint satisfaction problem where the constraint language is a finite-signature reduct of $(\Z;<)$. 
\end{definition}

\begin{example}\label{ex:F}
We give examples of reducts of $(\Z;<)$; 
the relations from those examples will re-appear 
in later sections.

\begin{enumerate}
\item $(\Z; \suc^p)$, where $\suc^p = \{(x,y)\in\Z^2 \mid y=x+p\}$ for $p\in\Z$,
\item $(\Z;\Diff{S})$, where $\Diff{S} := \{(x,y) : x,y \in \Z, y-x \in S\}$ for a finite set $S \subset Z$. 
\item $(\Z;\Diff{\{2\}}, \{(x,y) : |x-y| \leq 2\})$. 
\item $(\Z;F)$ where $F$ is the 4-ary relation $\{(x,y,u,v) : y=x+1 \Leftrightarrow v=u+1\}$. 
\item $(\Z;\neq,\Dist{i})$ where $\Dist{i} := \{(x,y) : |x - y| = i\}$.  
\item $(\Z; \{(x,y,z)\in\Z^3 \mid z\leq \max(x,y)+p\})$ for any $p\in\Z$,
\end{enumerate}
\end{example}

Examples 4.\ and 5.\ above are distance CSPs which have unbounded Gaifman degree (see Section~\ref{sect:degrees}), so they do not fall into the scope of~\cite{BodDalMarMotPin}.
If $p\neq 0$, the last example is a discrete temporal CSP that is not a temporal CSP and does not fall into the scope of~\cite{tcsps-journal}.
For a subset $S$ of the domain of a structure $\Gamma$, we write $\Gamma[S]$ for the structure induced on $S$ by $\Gamma$. 
It is easy to see that $(\Z;<)$ admits quantifier elimination in the language consisting of the binary relations $y<x+c$ and $y\leq x+c$ for $c\in\Z$.
We will call a formula \emph{standardized} if it is a quantifier-free formula in conjunctive normal form.
If $\phi$ is a first-order formula using literals of the form $x<y+c$ or $x\leq y+c$, we will allow ourselves to talk about $\phi$
as if it were a first-order formula in the language $\{<\}$, when this does not cause any confusion.
The following is easy to see. 

\begin{proposition}\label{prop:inNP}
All discrete temporal CSPs are in NP. 
\end{proposition}
\begin{proof}
Let $q$ be the size of the biggest integer that appears in the standardized formulas that define the relations in $\Gamma$ over $(\Z;<)$;
that is, for any atomic formula $x < y+k$ or $x \leq y+k$ in those formulas, $k \in \Z$, we have $|k| \leq q$. 
For an instance $\Phi$ of $\Csp(\Gamma)$ with $n$ variables, it is clear that $\Gamma \models \Phi$ if and only if $\Phi$ is true on  $\Gamma[\{1,\dots,(q+1)n\}]$.
We may guess a satisfying assignment of values from $\{1,\dots,(q+1)n\}$ to the variables of $\Phi$, 
and verify in polynomial time that all the constraints are satisfied. 
\end{proof}

\section{The Algebraic Approach}
\label{sect:algebra}
The starting point of the universal algebraic
approach to analyze the complexity of CSPs
is the observation that when a relation $R$ can be
defined by a  primitive positive formula over $\Gamma$,
then $\Csp(\Gamma)$ allows to simulate
the `richer' problem $\Csp(\Delta)$ where
$\Delta = (\Gamma,R)$ has been obtained from $\Gamma$ by adding $R$ as another relation. The proof of this fact given by Jeavons, Cohen, and Gyssens~\cite{JeavonsClosure} works for all structures $\Gamma$ over finite or over infinite domains. Since we will use this fact very frequently, we will not explicitly refer back to it from now on. 

\emph{Polymorphisms} are an important tool to study the question which
relations are primitive positive definable in $\Gamma$. We say that a function $f \colon D^n \to D$ \emph{preserves} a relation $R \subseteq D^m$ if for all $t_1,\dots,t_n \in R$
the tuple
$f(t_1,\dots,t_k)$ obtained by applying $f$ componentwise to the tuples $t_1,\dots,t_k$ is also in $R$; otherwise, $f$ \emph{violates} $R$. 
A \emph{polymorphism} of a relational structure $\Gamma$ with domain $D$ is a function 
from $D^n$ to $D$, for some finite $n$, 
which preserves \emph{all} relations of $\Gamma$. We write 
$\Pol(\Gamma)$ for the set of all polymorphisms of $\Gamma$. 

We write $\mathcal O$ for $\bigcup_{k \in \mathbb N} (D^k \to D)$. 
A subset $\mathcal F$ of $\mathcal O$ \emph{generates} $f \in \mathcal O$
if $f$ can be obtained from projections and functions in $\mathcal F$ by composition. 
Note that every function generated by polymorphisms of $\Gamma$
is again a polymorphism. 
We will need the fact that the set of all polymorphisms of $\Gamma$
 is furthermore \emph{locally closed}, 
 that is, when $f \in D^k \to D$ is such that for
all finite $S \subseteq D^k$ there exists
an $e \in \Pol(\Gamma)$ such that 
$e(x) = f(x)$ for all $x \in S$, then
$f$ is also a polymorphism of $\Gamma$. 
A subset $\mathcal F$ of $\mathcal O$ \emph{locally generates} $f\in\mathcal O$
if for every finite subset $S$ of $D$, there exists a function $g$
that is generated by $\mathcal F$ and such that the restrictions of $g$ and $f$
to $S$ coincide.

It is clear that a polymorphism
of a structure $\Gamma$ also preserves all relations that are primitive positive definable in $\Gamma$; this holds for arbitrary finite and infinite structures $\Gamma$. 
If $\Gamma$ is finite~\cite{BoKaKoRo,Geiger} or $\omega$-categorical~\cite{BodirskyNesetrilJLC},
then a relation is preserved by all polymorphisms \emph{if and only if} it is primitively positively definable
in $\Gamma$.

The structures that we consider in this paper will not be $\omega$-categorical;
however, following the philosophy in~\cite{BodHilsMartin-Journal}, one can refine these universal-algebraic methods to apply them also in our situation.
We will describe these refinements in the rest of this section.

The \emph{(first-order) theory} of a structure $\Gamma$, denoted by $\Th(\Gamma)$, is the set of all first-order sentences that are true in $\Gamma$.
We define some notation to conveniently work with 
models of $\Th(\Gamma)$ and their reducts. 

\begin{definition}[$\kappa.\Z$]
Let $\kappa$ be a linearly ordered set. We write $\kappa.\Z$ for $\kappa$ copies of $\Z$ indexed by the elements of $\kappa$; formally, $\kappa.\Z$ is the set $\{(a,z) : a \in \kappa, z \in \Z\}$. Then
$(\kappa.\Z;<)$ is the structure where $<$ denotes the lexicographic order on $\kappa.\Z$.
\end{definition}

It is well-known and easy to see that the models of 
$\Th(\Z; <)$ are precisely the structures isomorphic to $(\kappa.\Z;<)$, for some linear order $\kappa$. 
When $k \in \Z$ and $u = (a,z) \in \kappa.\Z$, we write $u + k$ for $(a,z+k)$. 

\begin{definition}[$\kappa.\Gamma$]\label{def:kappaGamma}
Let $\Gamma$ be a reduct of $(\Z;<)$ with signature $\tau$. Then $\kappa.\Gamma$ denotes
the `corresponding' reduct 
of $(\kappa.\Z;<)$ with signature $\tau$.
Formally, when $R \in \tau$ and 
$\phi_R$ is a formula
that defines $R^\Gamma$, then $R^{\kappa.\Gamma}$ is the relation defined by $\phi_R$ 
over $(\kappa.\Z;<)$. 
\end{definition}

In the following, we identify $\Z$ with the copy of $\Z$ induced by $0.\Z$ in $\Q.\Z$. That is, we view
$(\Z;<)$ as a substructure of $(\Q.\Z;<)$, and consequently $\Gamma$ 
as a substructure of $\Q.\Gamma$ for each reduct $\Gamma$ of $(\Z;<)$.

A \emph{type} of a structure $\Delta$
is a set $p$ of formulas with one free variable $x$ such that $p \cup \Th(\Delta)$ is satisfiable (that is, $\{\phi(c) : \phi \in p\} \cup \Th(\Delta)$, for a new constant symbol $c$, has a model). 
 A countable $\tau$-structure $\Gamma$ is \emph{saturated}
if for all choices of finitely many
constants $c_1,\dots,c_n$ for elements of $\Gamma$, 
and every unary type $p$ of $(\Gamma,c_1,\dots,c_n)$, 
there exists an element $d$ of $\Gamma$ such that $(\Gamma,c_1,\dots,c_n) \models \phi(d)$ for all $\phi \in p$. 
When $\Gamma$ and $\Delta$ are two countable saturated structures with the same first-order theory, then $\Gamma$ and $\Delta$ are isomorphic~\cite{Hodges}. 
Note that $(\Q.\Z;<)$ is saturated.
More generally, $\Q.\Gamma$ is saturated for every reduct $\Gamma$ of $(\Z;<)$. 

We define the function
$- \colon (\kappa.\Z)^2 \to (\Z \cup \{\pm\infty\})$
for $x,y \in \kappa.\Z$ by
\begin{align*}
x-y :=k \in \Z \quad & \text{ if } x = y+k \\
x-y := +\infty \quad & \text{ if } x \text{ and } y \text{ are in different copies and } x>y\\
x-y := -\infty \quad & \text{ otherwise.} 
\end{align*}

When $\Gamma$ and $\Delta$ are two structures with the same relational signature $\tau$, then a \emph{homomorphism} from $\Gamma$ to $\Delta$ is a function from the domain of $\Gamma$ to the domain of $\Delta$ such that for every $R \in \tau$ of arity $k$ we have
$R^\Gamma(u_1,\dots,u_k) \Rightarrow R^\Delta(f(u_1),\dots,f(u_k))$. 
It is straightforward to see that if there
is a homomorphism from $\Gamma$ to $\Delta$, and vice versa,
then $\Csp(\Gamma)$ and $\Csp(\Delta)$ are the same computational problem.

\begin{lemma}[See Lemma~2.1 in~\cite{BodHilsMartin-Journal}]\label{lem:sat}
Let $\Gamma$ be a countable saturated structure,
let $\Delta$ be countable, let $d_1,\dots,d_k$ be elements of $\Delta$,
and let $c_1,\dots,c_k$ be elements of $\Gamma$. 
Suppose that for all primitive positive formulas $\phi$ such that $\Delta \models \phi(d_1,\dots,d_n)$ we have 
$\Gamma \models \phi(c_1,\dots,c_k)$. 
Then there exists a homomorphism from $\Delta$ to $\Gamma$ that maps $d_i$ to $c_i$ for all $i \leq k$. 
\end{lemma}

An \emph{endomorphism} is a
unary polymorphism. To classify the computational complexity of the CSP for all reducts of a structure $\Gamma$, 
it often turns out to be important to 
study the possible endomorphisms of those reducts first, 
before studying the polymorphisms,
e.g. for the reducts of $({\mathbb Q};<)$ in
~\cite{tcsps-journal} and the reducts of the countably infinite random graph in~\cite{BodPin-Schaefer}. 

We are now in the position to state a general result, Theorem~\ref{thm:galois}, whose proof might explain the importance of saturated models for the universal-algebraic approach. 
When $\Gamma$ is a structure, then the \emph{orbit} of a $k$-tuple $(a_1,\dots,a_k)$ of elements of $\Gamma$
is the set $\{(\alpha(a_1),\dots,\alpha(a_n)) \mid \alpha \in \Aut(\Gamma)\}$. 

\begin{theorem}\label{thm:galois}
Let $\Gamma$ be a countable saturated structure, let $\Delta$ be a reduct of $\Gamma$, and $R$ a relation with a first-order definition in $\Gamma$. Then 
\begin{itemize}
\item $R$ has a first-order definition in $\Delta$ if and only if $R$ is preserved by the automorphisms of $\Delta$;
\item $R$ has an existential positive definition in $\Delta$ if and only if $R$ is preserved by all the endomorphisms of $\Delta$;
\item if $R$ consists of $n$ orbits of $k$-tuples in $\Delta$, then $R$ has a primitive positive definition in $\Delta$ if and only if $R$ is preserved by all polymorphisms of $\Delta$ of arity $n$.
\end{itemize}
\end{theorem}
\begin{proof}
Suppose that $R$ is $k$-ary, and let $\phi$ be the first-order definition of $R$ in $\Gamma$.
It is well-known and straightforward to prove that first-order formulas
are preserved by automorphisms of $\Delta$, that existential positive formulas are preserved by endomorphisms of $\Delta$,
and that primitive positive formulas are preserved by polymorphisms of $\Delta$. 

Suppose first that $R$ is preserved by all automorphisms of
$\Delta$. Let $\phi$ be a first-order definition of $R$ in $\Gamma$.
Let $\Psi$ be the set of all first-order formulas in the language of $\Delta$ that are consequences of $R$.
Formally,
\[\Psi = \{\psi(x_1,\dots,x_k) \mid \forall (a_1,\dots,a_k)\in R, \Delta\models\psi(a_1,\dots,a_k)\}.\]
We prove that if a tuple $\tuple a$ satisfies every formula in $\Psi$ then $\tuple a$ is in $R$.
Let $\tuple a$ be such a tuple.
Let $p$ be the type of $\tuple a$ in $\Delta$.
By replacing in $p$ every relation symbol of the signature of $\Delta$
by a first-order definition of the corresponding relation in $\Gamma$, we obtain a set $q$ of formulas in the language of $\Gamma$.
If we can find some tuple $\tuple b$ that satisfies $\{\phi\}\cup q$ in $\Gamma$, then we are done. Indeed, we have that $\tuple b$ is in $R$,
and $\tuple b$ has the same type as $\tuple a$ in $\Delta$. By saturation of $\Delta$ and a back-and-forth argument, $\tuple b$ and $\tuple a$ are in the same orbit under $\Aut(\Delta)$
so that $\tuple a$ is in $R$ as well.
So let us assume that $\{\phi\}\cup q$ is not satisfiable in $\Gamma$. 
By saturation of $\Gamma$, there is some formula $\psi$ in $q$ such that $\{\psi,\phi\}$ is not satisfiable in $\Gamma$,
i.e., we have that $\Gamma\models \forall x_1,\dots,x_k (\phi(x_1,\dots,x_k)\implies\neg\psi(x_1,\dots,x_k))$.
The formula $\psi$ corresponds to a formula $\theta$ in the language of $\Delta$ by construction.
We therefore obtain that $\neg\theta\in\Psi$, so $\neg\theta\in p$. But $\theta\in p$, a contradiction.

Suppose now that $R$ is preserved by all endomorphisms of $\Delta$.
In particular $R$ is preserved by all the automorphisms of $\Delta$, so that there exists a first-order definition $\phi$ of $R$ in $\Delta$.
Let $\Psi$ be the set of all universal negative consequences of $\phi$ in $\Delta$.
As above, we aim to prove that if $\tuple a$ satisfies all the formulas in $\Psi$, then $\tuple a$ is in $R$.
Let $\tuple a$ be such a tuple, and let now $p$ be the ep-type of $\tuple a$, that is, the set of all the existential positive formulas $\psi$
such that $\Delta\models\psi(\tuple a)$.
If $p\cup\{\phi\}$ is satisfiable in $\Delta$, then we are done: there exists a tuple $\tuple b\in R$ that has the same ep-type as $\tuple a$.
By saturation of $\Delta$, we can produce an endomorphism of $\Delta$ that maps $\tuple b$ to $\tuple a$, so that $\tuple a$ is in $R$.
Otherwise, there exists a single formula $\theta\in p$ such that $\Delta\models\forall x_1,\dots,x_k (\phi(x_1,\dots,x_k)\Rightarrow\neg\theta(x_1,\dots,x_k))$.
Therefore $\neg\theta$ is in $\Psi$, so that $\tuple a$ must satisfy $\neg\theta$, contradicting the fact that $\tuple a$ already satisfies $\theta$.

Finally, suppose that $R$ consists of $n$ orbits of $k$-tuples in $\Gamma$, 
and that $R$ is preserved by all polymorphisms of $\Delta$ of arity $n$. Let $\Psi$ be the set of all primitive positive formulas with free variables
$x_1,\dots,x_n$ that hold for all tuples in $R$ in $\Delta$, 
and let $(a_1,\dots,a_k)$ be
a tuple that satisfies $\Psi$. 
Pick representatives $(b_1^1,\dots,b_k^1),\dots,(b_1^n,\dots,b_k^n)$ for all the orbits of $k$-tuples in $\Gamma$ that lie in $R$. 
Note that every primitive positive formula that holds on 
$(b^1_1,\dots,b^n_1),\dots,(b^1_k,\dots,b^n_k)$ in $\Delta^n$ also holds
on $(a_1,\dots,a_k)$ in $\Delta$. 
By Lemma~\ref{lem:sat} and saturation of $\Delta$, 
there exists a homomorphism
from $\Delta^n$ to $\Delta$ that maps $(b^1_i,\dots,b^n_i)$ 
to $a_i$ for all $i \in \{1,\dots,k\}$. This map is a polymorphism of $\Delta$, and since $R$ is preserved by
polymorphisms, $(a_1,\dots,a_n) \in R$. 
Similarly as before, a compactness argument for 
first-order logic over $\Gamma$ shows that
$\Psi$ is equivalent to a single primitive positive sentence that is equivalent to $\phi$. 
\end{proof}

\section{Statement of Results}\label{sect:results}
In this section, we describe the border between NP-complete discrete temporal CSPs and polynomial-time tractable discrete temporal CSPs, modulo the Feder-Vardi dichotomy conjecture.

	\begin{definition}
	Let $d$ be a positive integer. 
	The \emph{d-modular max}, $\max_d\colon\Z^2\to\Z$,
	is defined by $\max_d(x,y):=\max(x,y)$ if $x=y\mod d$ and $\max_d(x,y):=x$ otherwise.
	The \emph{d-modular min} is defined analogously.  
	\end{definition}
	
	Note that these two operations are not commutative when $d>1$.
    Examples of relations which are preserved by $\max$ and which are definable over $(\Z;<)$ are the relations appearing in the last item of Example~\ref{ex:F}.
    Another example of a relation which is preserved by $\max_d$ is the ternary relation containing the triples of the form
    $(a+d,a,a),(a+d,a+d,a),(a,a+d,a)$
    for all $a\in\Z$. Note that for a fixed $d$, the relation of the latter type is preserved by $\max_d$ but not by $\max_{d'}$ for any other $d'$.

\begin{theorem}\label{thm:main}
Let $\Gamma$ be a reduct of $(\Z;<)$ with finite signature.
Then there exists a structure $\Delta$ such that $\Csp(\Delta)$ equals
$\Csp(\Gamma)$ and one of the following cases applies. 
\begin{enumerate}
\item $\Delta$ has a finite domain, and the CSP for $\Gamma$ is conjectured to be in Ptime or NP-complete~\cite{FederVardi}.
\item $\Delta$ is a reduct of $(\Q;<)$, and the complexity of $\Csp(\Delta)$ has been classified in~\cite{tcsps-journal}. 
\item $\Delta$ is a reduct of $(\Z;<)$ and preserved by a modular max or modular min.
In this case, $\Csp(\Gamma)$ is in Ptime.
\item $\Delta$ is a reduct of $(\Z;\suc)$ that is preserved by a binary injective function preserving $\suc$. 
In this case, $\Csp(\Gamma)$ is in Ptime. 
\item $\Csp(\Gamma)$ is NP-complete. 
\end{enumerate}
\end{theorem}

As an illustration of the algorithmic consequences of our main result, we give
examples of computational problems that can be formulated as discrete temporal CSPs and that are in Ptime (the corresponding constraint languages are preserved 
by the $d$-modular maximum function). Fix positive integers $d$ and $K$. 

\medskip \noindent
{\bf Input:} $A \in \{0,1\}^{m \times n}$,
$B \in \{0,1\}^{m' \times n}$, and $c \in \{-K,\dots,K\}^m$ such that in any row of $A$ and $B$ there is precisely one $1$ and one $-1$. \\
{\bf Question:} Is there a vector $x\in\Z^n$ such that $Ax \leq c$ and $Bx=0 \bmod d$? 

\medskip
This problem can be seen as $\Csp(\Z; B_{-K},\dots,B_K, \Mod_d)$ 
where $B_i := \{(x,y) \in\Z^2 \mid y-x\leq i\}$ and $\Mod_d := \{(x,y)\in\Z^2 \mid x=y\bmod d\}$.

\medskip


\medskip 
\noindent{\bf Input:} 
A directed graph $(V,E)$
and a weight function 
$f \colon E \to \{1,d\}$
such that $f(x,y)=d$ implies that 
there is exactly one 
other edge of weight $d$ leaving $x$. 
\\ 
{\bf Question:} 
Can we remove at most one outgoing
edge $(x,y)$ of weight $d$ per vertex $x$ so that in the resulting graph 
every cycle has algebraic weighted sum 0? 

\medskip
This problem 
is essentially 
$\Csp(\Z; \suc,R)$ where
$$R = \{(x,y_1,y_2) \in \Z^3 \mid 
y_1-x, y_2-x \in \{0,d\}  \text{ and } y_1 = x + d \vee y_2 = x + d\}.$$
Indeed, given an instance as above, 
we can add a variable for each vertex of the graph, a constraint $R(x,y_1,y_2)$
for each pair of edges $(x,y_1),(x,y_2) \in E$ of weight $d$, and a $\suc$ constraint for all edges of weight one. The original instance is satisfiable iff the created instance has a solution. 
In this case, we can find the edges that should be removed  
from the equalities $y_i=x+d$ that are not satisfied in a solution to the new instance.

\section{Definability of Successor and Order}
\label{sect:succ}
The goal of this section is a proof that the CSPs for reducts of $(\Z;<)$ fall into four classes.
This will allow us to focus in later sections on reducts of $(\Z;<)$ where $\suc$ is pp-definable.

\begin{theorem}\label{thm:class}
Let $\Gamma$ be a reduct of $(\Z;<)$ with finite signature. 
Then $\Csp(\Gamma)$ equals $\Csp(\Delta)$
where $\Delta$ is one of the following:
\begin{enumerate}
\item a finite structure;
\item a reduct of $(\Q;<)$;
\item a reduct of $(\Z;<)$ where $\Dist{k}$ 
is pp-definable for all $k \geq 1$;
\item a reduct of $(\Z;<)$ where $\suc$ is pp-definable.
\end{enumerate}
\end{theorem}

The proof of this result requires some effort and
spreads over the following subsections.
Before we go into this, we explain the significance
of the four classes for the CSP. 

It is easy to see that there exists a structure
$\Delta$ with a finite domain such that $\Csp(\Gamma)$ equals
$\Csp(\Delta)$ if and only if $\Gamma$ has an endomorphism with finite range. So we will assume in the following that this is not the case. 
The CSPs for reducts of $(\Q;<)$ have been studied in~\cite{tcsps-journal}; they are either in Ptime or NP-complete. Hence, 
we are also done if there exists a reduct $\Delta$ of 
$(\Q;<)$ such that $\Csp(\Delta) = \Csp(\Gamma)$.
Several equivalent characterizations of those reducts $\Gamma$ will be given in Section~\ref{sect:petrus}. 
This is essential for proving Theorem~\ref{thm:class}. 
When $\Gamma$ is a reduct of $(\Z;<)$ where for all $k \geq 1$ the relation $\Dist{k}$ is pp-definable,
then $\Csp(\Gamma)$ is NP-complete; this is a consequence of the following proposition from~\cite{BodDalMarMotPin}.

\begin{proposition}[Proposition~27 in~\cite{BodDalMarMotPin}]
\label{prop:26}
Suppose that the relations $\Dist1$ and $\Dist5$ are pp-definable in $\Gamma$. Then $\Csp(\Gamma)$ is NP-hard.
\end{proposition}

The previous paragraphs explain why  Theorem~\ref{thm:class} indeed  reduces the complexity classification 
of CSPs for finite-signature reducts $\Gamma$ of $(\Z;<)$ to the case where $\suc$ is pp-definable in $\Gamma$. 

\subsection{Degrees}\label{sect:degrees}
We consider three notions of \emph{degree} for relations $R$ that are first-order definable in $(\Z;<)$:
\begin{itemize}
\item For $x \in \Z$, we consider the number of $y \in \Z$ that appear together with $x$ in a tuple from $R$;
this number is the same for all $x \in \Z$, and called
the \emph{Gaifman-degree} of $R$ (it is the degree of the Gaifman graph of $({\mathbb Z};R)$). 
\item The \emph{distance degree} of $R$ is the supremum of $d$ such that there are $x,y \in \mathbb Z$ that occur together in a tuple of $R$ and $|x-y|=d$. 
\item The \emph{quantifier-elimination-degree (qe-degree)} of $R$ is the minimal $q$ so that there is a quantifier-free definition $\phi$ of $R$,
such that for every literal $x<y+k$ or $x\leq y+k$ in $\phi$,
we have $|k|\leq q$.
\end{itemize}
The degree of a reduct of $(\Z;<)$ is the supremum of the degrees of its relations, for any of the three notions of degree.
The paper \cite{BodDalMarMotPin} considered reducts of  $(\mathbb{Z};\suc)$ with finite Gaifman-degree.
Note that the Gaifman-degree is finite if and only if the distance degree is finite.
In this paper, qe-degree will play the central role, as any reduct of $(\mathbb{Z};<)$ with finite relational signature clearly has finite qe-degree.



\subsection{Compactness}

In this section we present some results, based 
on applications of K\"onig's tree lemma,
that show how properties about finite substructures
of finite-signature reducts $\Gamma$ of $(\Z;<)$ correspond to the existence of certain homomorphisms from $\Gamma$ to $\Q.\Gamma$.

Let $(\kappa.\Z;<)$ be a model
of $\Th(\Z;<)$, let $S$ be any set, let $s \in \N$, 
and $f \colon S \to \kappa.\Z$. 
We say that $x,y \in S$ are \emph{$(f,s)$-connected} if there is a sequence $x = u_1,\ldots,u_k = y \in S$ so that $0 \leq |f(u_i)-f(u_{i+1})| \leq s$ 
for all $i \in \{1,\dots,k-1\}$.
Note that this notion of connectivity defines an equivalence relation on $S$.
We define an equivalence relation $\sim_s$ on
functions $f,g \colon S \to \kappa.\Z$ as follows: $f \sim_s g$ when 
the following conditions are met:
\begin{itemize}
\item $x,y \in S$ are $(f,s)$-connected
if and only if they are $(g,s)$-connected,
\item if $x,y \in S$ are $(f,s)$-connected (and therefore $(g,s)$-connected) then $f(x)-f(y)=g(x)-g(y)$,
\item if $x,y\in S$ are not $(f,s)$-connected then $f(x)<f(y) \Leftrightarrow g(x)<g(y)$.
\end{itemize}
In other words, $f\sim_s g$ iff the equivalence relations defined by $(f,s)$-connectivity and $(g,s)$-connectivity have the same equivalence classes,
are such that within each equivalence class the pairwise distances are the same, and the order of the equivalence classes are the same.

\begin{lemma}[Substitution Lemma]\label{lem:substitution}
Let $\Gamma$ be a reduct of $(\Z;<)$
with qe-degree $q$.
Let $S \subseteq \kappa_1.\Z$, 
and $f,g \colon S \rightarrow \kappa_2.\Z$ 
be such that $f \sim_q g$. 
Then $f$ is a homomorphism from 
$(\kappa_1.\Gamma)[S]$ 
to $\kappa_2.\Gamma$
if and only if $g$ is such a homomorphism. 
\end{lemma}
\begin{proof}
Let $\phi$ be a quantifier-free definition of a $k$-ary relation from $\Gamma$, and let $\bar c = (c_1,\dots,c_k) \in (\kappa_1.\Z)^k$ be such that $(\kappa_1.\Z;<) \models \phi(\bar c)$. 
Since the qe-degree of
$\Gamma$ is $q$, we can assume that all atomic formulas of $\phi$ are of the form
$y \leq x+i$ or of the form $y < x+i$, for $i$ an integer such that $|i|\leq q$.
Suppose that $f$ is a homomorphism from $(\kappa_1.\Gamma)[S]$ to  $\kappa_2.\Gamma$. 
Since $f$ preserves $\phi$, we have
$(\kappa_2.\Z;<) \models \phi(f(\bar c))$. 
Now, let $i$ be such that $|i|\leq q$. 
Note that if $x,y \in \{c_1,\dots,c_n\}$ are $(f,q)$-connected,
then $f(x) \leq f(y)+i$ if and only if 
$g(x) \leq g(y)+i$, since $f \sim_q g$.
If $x,y$ are not $(f,q)$-connected we have $f(x)<f(y)\Leftrightarrow g(x)<g(y)$.
If $f(x)\leq f(y)+i$ holds with $i\geq 0$ then $f(x)< f(y)$
(otherwise $|f(x)-f(y)|\leq q$, contradicting the fact that $x$ and $y$
are not $(f,q)$-connected),
so $g(x)< g(y)$ holds and $g(x)\leq g(y)+i$.
If $i<0$ we obtain $g(x)<g(y)$, and since $x,y$ are not $(g,q)$-connected,
we have $g(x)<g(y)-q \leq g(y)+i$.
The proof for literals of the form $x<y+i$ is similar.

It follows that $(\kappa_2.\Z;<) \models \phi(g(\bar c))$, and hence $g$ is a homomorphism from $(\kappa_1.\Gamma)[S]$ to $\kappa_2.\Gamma$ as well. 
\end{proof}

\begin{proposition}[Compactness]\label{prop:compactness}
   Let $S$ be a subset of $\Q.\Z$ and let $(a_i)_{i\in\N}$ be an enumeration of $S$.
   Let $s\geq 0$.
   Let $(F_i)_{i\in\N}$ be a sequence of sets such that:
   \begin{enumerate}
    \item for all $i$, $F_i$ is the $\sim_s$-equivalence class of some function $\{a_0,\dots,a_i\}\to\Q.\Z$,
    \item if $g\in F_j$ and $i<j$, then $g|_{\{a_0,\dots,a_i\}}\in F_i$.
   \end{enumerate}
   Then there exists a function $h\colon S\to\Q.\Z$ such that for all $i$,
   we have $h|_{\{a_0,\dots,a_i\}}\in F_i$.
\end{proposition}
\begin{proof}
    We define the function $h$ by induction.
    We require that at each step, the function $$h_i\colon\{a_0,\dots,a_i\}\to\Q.\Z$$ that we define is in $F_i$ and is such that whenever $a,b\in\{a_0,\dots,a_i\}$
    are not $(g,s)$-connected for \emph{any} function $g$ in some $F_j$, then $h_i(a)-h_i(b)=\infty$. Conversely if $a,b\in\{a_0,\dots,a_i\}$
    are $(g,s)$-connected for \emph{some} function $g\in F_j$, then $h_i(a)-h_i(b)=g(a)-g(b)$.
    For $i=0$, let $h_0$ be any function in $F_0$.
    Suppose now that $h_i$ has been defined, and let $h_{i+1}(a_j)=h_i(a_j)$ for $j\in\{0,\dots,i\}$.
    Let $g\in F_j$ be such that for every pair $a_k,a_l\in\{a_0,\dots,a_{i+1}\}$,
    if $a_k,a_l$ are $(g',s)$-connected for some $g'\in F_{j'}$ then they are actually $(g,s)$-connected:
    such a function exists, by taking $j$ sufficiently large that $\{a_0,\dots,a_j\}$
    contains all the elements that witness that $a_k,a_l$ are $(g',s)$-connected for some $g'$.
    From the induction hypothesis and the properties 1.\ and 2., we know that $h_{i}\sim_s g|_{\{a_0,\dots,a_i\}}$.
    Define $h_{i+1}(a_{i+1})$ as follows:
    \begin{itemize}
        \item suppose that there exists $k$ such that $a_{i+1}$ and $a_k$ are $(g,s)$-connected.
            Define $h_{i+1}(a_{i+1}) = h_i(a_k) - g(a_k) + g(a_{i+1})$.
        \item otherwise consider the sets
        $$P=\{p\in\Q \mid \exists k\in \{0,\dots,i\} : g(a_{k}) < g(a_{i+1}) \text{ and } h_i(a_k)\in p.\Z\}$$
        and
        $$Q=\{q\in \Q\mid \exists k\in \{0,\dots,i\} : g(a_{i+1})<g(a_{k}) \text{ and } h_i(a_k)\in q.\Z\}.$$
        We have $P<Q$. Indeed, let $p\in P,q\in Q$, and let $k,l\in\{0,\dots,i\}$ be such that $h_i(a_k)\in p.\Z$ with $g(a_k)<g(a_{i+1})$
        and $h_i(a_{l})\in q.\Z$ with $g(a_{i+1})<g(a_l)$.
        Since $a_{i+1}$ is not $(g,s)$-connected to some element of $\{a_0,\dots,a_i\}$,
        we have that $a_k$ and $a_l$ are not $(g,s)$-connected.
        By construction, we therefore have that $h_i(a_k)-h_i(a_l)=\infty$.
        Since $a_k$ and $a_l$ are not $(g,s)$-connected and since $g(a_k)<g(a_l)$, we have that $h_i(a_k)<h_i(a_l)$.
        It follows that $p<q$.
        Thus, there exists a rational $p$ such that $P<p<Q$. Define $h_{i+1}(a_{i+1}) = (p,0)$.
    \end{itemize}
    We now prove that the induction hypothesis remains true for $h_{i+1}$. 
    We claim that $h_{i+1}\sim_s g|_{\{a_0,\dots,a_{i+1}\}}$.
    Remember that we already now that $h_{i}\sim_s g|_{\{a_0,\dots,a_i\}}$ since $h_i\in F_i$ by induction and $g\in F_j$ for $j>i$.
    Let $a_j\in\{a_0,\dots,a_i\}$. If $h_{i+1}(a_{i+1})$ is at finite distance from $h_{i+1}(a_j)$,
    then by definition $a_j,a_{i+1}$ are $(g,s)$-connected. Let $k\in\{0,\dots, i\}$ be the index used in the definition of $h_{i+1}$.
    We then have
    \begin{align*}
    h_{i+1}(a_{i+1}) - h_{i+1}(a_j) &= h_i(a_k)-g(a_k)+g(a_{i+1})-h_i(a_j)\\
                        &= g(a_k)-g(a_j) - g(a_k) + g(a_{i+1})\\
                        &= g(a_{i+1})-g(a_j).
    \end{align*}
    If $h_{i+1}(a_{i+1})$ and $h(a_j)$ are at infinite distance, then $a_{i+1},a_j$ are not $(g,s)$-connected.
    The order induced on $a_{i+1}$ and $a_j$ by $h_{i+1}$ is then the same as the order induced by $g$,
    from the construction of $P$ and $Q$. It follows that $h_{i+1}\sim_s g|_{\{a_0,\dots,a_{i+1}\}}$.
    Moreover $h_{i+1}$ indeed separates integers that are never $(g,s)$-connected for any $g\in F_j$.
    Finally, if $g'\in F_{j'}$ is such that $a,b$ are $(g',s)$-connected then $a$ and $b$ are also $(g,s)$-connected
    and $g'(a)-g'(b)=g(a)-g(b)$, so that the choice of $g$ in our construction is irrelevant.
    This proves that $h_{i+1}$ satisfies the induction hypothesis.

    It remains now to take $h=\bigcup_{i\geq 0} h_i$, which satisfies the conclusion of the statement.
\end{proof}

Compactness and substitution will be applied frequently; one application is in the proof of the following proposition. Note that this makes essential use of the saturated model.

\begin{proposition}
\label{prop:koenig}
Let $\Gamma$ be a finite-signature 
reduct of $(\Z;<)$. 
Then for all $a_1, a_2 \in \Z$ either
\begin{itemize}
\item there is an $r \geq 0$ and a finite $S \subseteq \Z$ that contains 
$\{a_1,a_2\}$ such that for all homomorphisms $f$ from $\Gamma[S]$ to $\Gamma$ we have $|f(a_1)-f(a_2)| \leq r$, or
\item  there is a homomorphism $h$ from
$\Gamma$ to $\Q.\Gamma$
 such that $h(a_1)-h(a_2) = \infty$. 
\end{itemize}
\end{proposition}

\begin{proof}
Let  $a_1,a_2 \in \Z$ be arbitrary. 
Suppose that for all $r \geq 0$ 
and all finite $S \subset \Z$ 
containing $\{a_1,a_2\}$ there is a homomorphism $f$ from $\Gamma[S]$ 
to $\Gamma$ 
such that $|f(a_1)-f(a_2)| > r$. 
We will describe how to construct the desired homomorphism $h$.

Let $a_1,a_2,a_3,\dots$ be an enumeration of $\Z$, and let $q$ be the qe-degree of $\Gamma$. 
Consider the following infinite tree $\mathcal{T}$ whose vertices lie on levels $1,2,\dots$
The vertices at the $n$-th level are the $\sim_q$-equivalence classes of homomorphisms $f$ from 
$\Gamma[\{a_1,\dots,a_{n+1}\}] \to \Q.\Gamma$ that 
satisfy $|f(a_1)-f(a_2)| > qn$ where $q$ is the qe-degree of $\Gamma$. 
We have an arc in $\mathcal T$ from an equivalence class $F$ on level $n$ to an equivalence class $G$
on level $n+1$ if there are $f \in F$, $g \in G$
such that $f$ is the restriction of $g$. 
By assumption, $\mathcal{T}$ has vertices on each level $n$. The tree
$\mathcal{T}$ has finitely many vertices on each level, since the number of $\sim_q$-equivalence classes of homomorphisms 
from $\Gamma[\{a_1,\dots,a_n\}] \to \Q.\Gamma$ is bounded by $n^2(q+2)$. 

It follows by K\"onig's lemma that there is an infinite branch $\mathcal{B}$ of $\mathcal{T}$.
By Proposition~\ref{prop:compactness} using the elements of $\mathcal B$ for the sequence $(F_i)_{i\in\N}$,
there exists a function $h\colon \Z\to\Q.\Z$ such that for every $i\in\N$, $h|_{\{a_1,\dots,a_i\}}$ is in the branch $\mathcal B$.
Therefore we have $h(a_1)-h(a_2)=\infty$ and moreover $h$ is a homomorphism $\Gamma\to\Q.\Gamma$ by Lemma~\ref{lem:substitution}.
\end{proof}

\begin{definition}
A mapping $h$ between models of the
first-order theory of $(\Z;<)$ is called 
\emph{isometric} if $|h(x)-h(y)|=|x-y|$ 
for all $x,y \in \Z$.
\end{definition}

The following proposition can be shown by 
straightforward modifications of the proof 
of Proposition~\ref{prop:koenig}. 

\begin{proposition}\label{prop:isometric}
Let $\Gamma$ be a finite-signature 
reduct of $(\Z;<)$. Then either
\begin{itemize}
\item for all $r$ there is a finite $S \subseteq \Z$ containing $\{0,r\}$ such that for all homomorphisms $f$ from $\Gamma[S]$ to $\Gamma$ we have $|f(0)-f(r)| = r$, or
\item there is a homomorphism $h$ from $\Gamma $ to $\Q.\Gamma$ which is not isometric. 
\end{itemize}
\end{proposition}

\subsection{Finite-range endomorphisms}

In this section we present a lemma that
gives a useful sufficient condition for
$\Gamma$ to have endomorphisms with finite range. We will need the following combinatorial definitions and lemmas about the integers. 

We say that $T \subseteq \Z$ \emph{contains
arbitrarily long intervals} when for all $m \in \N$ there exists $z \in \Z$ so that $[z,z+m] \subset T$. 
A sequence $u_1,\dots,u_r$ is called a
\emph{\mprog} if $1 \leq u_{i+1} - u_i \leq m$
for all $i < r$. 
We say that $T$ has \emph{arbitrarily long \mprog s} if for all $r \in \N$ the set $T$ contains a \mprog~$u_1,\dots,u_r$.
Clearly, 
if $\Z \setminus T$ 
does not have arbitrarily long 
intervals then there exists an $m \in \N$ so that $T$
has arbitrarily long \mprog s.

\begin{lemma}\label{lem:partition}
Let $T \subseteq \Z$ contain
arbitrarily long \mprog s,
and let $T = T_1 \cup \dots \cup T_k$ be a partition of $T$ into finitely many sets. Then there exists an $i \leq k$ and an $m' \in \N$ such that $T_i$ contains arbitrarily long 
$({\leq}m')$-progressions.
\end{lemma}

\begin{proof}
If there exists an $m' \in \N$ such that 
$T_1$ contains arbitrarily long $({\leq}m')$-progressions, then there is nothing to show. 
So suppose that this is not the case. 

We will show that $T' := T \setminus T_1$
contains arbitrarily long $({\leq}m)$-progressions; the statement then clearly follows by induction. 
Let $s \in \N$ be arbitrary. We want 
to find a \mprog~$u_1,\dots,u_s$ in $T'$. 
By the above assumption, 
$T_1$ does not contain arbitrarily 
long $({\leq}ms)$-progressions,
and hence there exists an $r$
such that $T_1$ does not contain 
a $({\leq}ms)$-progression of length $r$.

Since $T$ contains arbitrarily long \mprog s,
it contains in particular an \mprog~$\rho$ of length $msr$. Consider the first $s$ elements 
of $\rho$. If all those 
elements are in $T'$ we have found
the desired \mprog~of length $s$, and are done. So suppose otherwise; that is, at least one of those first $s$ elements must be from $T_1$. 
We apply the same argument to the
next $s$ elements of $\rho$,
and can again assume that at least one
of those elements must be from $T_1$.
Continuing like this, we find a subsequence
of $\rho$ of elements of $T_1$
which form a $({\leq}ms)$-progression.
The length of this subsequence is 
$msr/ms = r$. 
But this contradicts our assumption
that $T_1$ does not contain $({\leq}ms)$-progression of length $r$. 
\end{proof}

\begin{lemma}\label{lem:pigeon-hole}
Let $m \in \N$ and let $T \subseteq \Z$ be with arbitrarily long \mprog s. Then for all $S \subset \Z$ of cardinality $m+1$ 
there are $x_1,x_2 \in S$ and $y_1,y_2 \in T$ such that $x_1 - x_2 = y_1 - y_2$. 
\end{lemma}
\begin{proof}
Let $r$ be greater than $\max(S)-\min(S)$. 
Then there exists an \mprog \ $w_1,\dots,w_r$ in $T$.
Define $T_i := \{z - w_1 + \min(S) + i \; | \; z \in T\}$. 
Then $T_0 \cup \cdots \cup T_{m-1}$ includes
the entire interval $[\min(S),\max(S)]$.
By the pigeon-hole principle there is an $i$ such that
$|T_i \cap S| \geq 2$, which clearly implies the statement.
\end{proof}

\begin{lemma}\label{lem:koenig2}
Let $\Gamma$ be a finite-signature reduct 
of $(\Z;<)$, let $h$ be a homomorphism 
from $\Gamma \to \kappa.\Gamma$ for some $\kappa$, and let $S \subseteq \Z$ be finite. Let $z_0\in\kappa.\Z$. 
If $\Z \setminus h^{-1}(S) \cap \{z\in\kappa.\Z : z\geq z_0\}$ does not contain arbitrarily long intervals then $\Gamma$ has a finite-range endomorphism. 
\end{lemma}
\begin{proof}
Since $\Z \setminus h^{-1}(S)$ does not contain
arbitrarily long intervals, there exists an $m'$ such that $T := h^{-1}(S)$ contains arbitrarily long $({\leq}m')$-progressions. Suppose that 
$S = \{s_1,\dots,s_k\}$, 
and define $T_i := h^{-1}(s_i)$.
Then by Lemma~\ref{lem:partition}
there exists an $m \in \mathbb N$ and an $i \leq k$ such that 
$T_i$ contains arbitrarily long $({\leq}m)$-progressions. 

Our argument is based on K\"onig's tree lemma, involving a finitely branching infinite tree $\mathcal T$, but with subtle differences when compared to
the construction given in the proof of Proposition~\ref{prop:koenig}. Let
$a_1,a_2,\dots$ be an enumeration of $\Z$,
and let $q$ be the qe-degree of $\Gamma$.  The vertices of $\mathcal T$ on the $n$-th level are the $\sim_q$-equivalence classes of 
homomorphisms $g$ from $\Gamma[\{a_1,\dots,a_n\}]$ to $\Gamma$ such that  $|g(\{a_1,\dots,a_n\})| \leq m$.
Adjacency is defined by restriction,
and $\mathcal T$ is finitely branching,
 as in the proof of Proposition~\ref{prop:koenig}.

We show that $\mathcal T$ has vertices on all
levels $n$ by induction on $n$.
We in fact prove the stronger statement that for any finite set $X\subset\Z$, there exists a homomorphism $g\colon\Gamma[X]\to\Gamma$ whose range has size at most $m$.
For $|X|\leq m$, this is witnessed by the restriction of the identity function to $X$.
Now let $|X|=n+1, n\geq m$. By Lemma~\ref{lem:pigeon-hole}, there are $x_j,x_k\in X$ and $y_1,y_2\in T_i$ such that $x_j-x_k = y_1-y_2$.
We therefore have that $f\colon x\mapsto h(x-x_j+y_1)$ is a homomorphism $\Gamma[X]\to\kappa.\Gamma$ whose range has size at most $n$.
Indeed, we have $f(x_j)=h(y_1)=h(y_2) = h(x_k-x_j+y_1)=f(x_k)$.
Up to $\sim_q$-equivalence, we can replace $f$ by another homomorphism $f'\colon\Gamma[X]\to\Gamma$. Let now $g$ be given by the induction hypothesis,
with $X=\mathrm{im}(f')$. We then have that $g\circ f'$ is a homomorphism $\Gamma[X]\to\Gamma$ whose range has size at most $m$, and the claim is proved.

Hence, $\mathcal T$ has vertices on all levels,
and therefore an infinite branch $\mathcal B$ by K\"onig's lemma. 
By Proposition~\ref{prop:compactness} and Lemma~\ref{lem:substitution}, we have a homomorphism $e\colon\Gamma\to\Gamma$
whose range has cardinality at most $m$, concluding the proof.
\end{proof}

\subsection{Endomorphisms of the saturated model}
This section presents an important lemma to 
analyze the endomorphisms of the saturated model $\Q.\Gamma$ of finite-signature reducts $\Gamma$
of $({\mathbb Z};<)$. 
The following lemma is already interesting and useful if $h$ is an endomorphism of $\Gamma$;
however, in some situations we need that lemma 
for homomorphisms from $\Gamma$ to $\Q.\Gamma$; however, in this case we have to make the additional assumption that
the set $\{q\in\Q \mid \exists z\in\Z : h(z)\in q.\Z\}$ of copies that are touched by the image of $h$ is bounded.
This is usually not a very restrictive assumption, since there is an $e \in \End(\Q.\Z;<)$ which has this property,
and hence for any $h \colon \Z \to \Q.\Z$ the mapping $e \circ h$ also has this property. 

\begin{lemma}\label{lem:hom-extension}
Let $\Gamma$ be a finite-signature reduct of $(\Z;<)$ without finite-range endomorphisms,
and let $h$ be a homomorphism
from $\Gamma$ to $\Q.\Gamma$
such that the set $\{q\in\Q\mid \exists z\in\Z : h(z)\in q.\Z\}$ is bounded. 
Then there exists an $e \in \End(\Q.\Gamma)$ which extends $h$ 
such that for all $x,y \in \Q.\Z$ with $x-y = \infty$ we have $e(x)-e(y) = \infty$. 
\end{lemma}
\begin{proof}
As in the proof of Proposition~\ref{prop:koenig}, we build $e$ through an argument involving 
K\"onig's lemma, the quantifier-degree $q$ of $\Gamma$, and an infinite tree $\mathcal T$. Let $a_1,a_2,\dots$
be an enumeration of $\Q.\Z$. In the $n$-th
level of $\mathcal T$ we will consider $\sim_q$-classes of homomorphisms $f$ from
$\Q.\Gamma[\{a_1,\dots,a_n\}]$ to $\Q.\Gamma$ with the property that 
\begin{itemize}
\item for all $x,y \in \{a_1,\dots,a_n\}$ with $x-y = \infty$ we have $f(x)-f(y) = \infty$, 
and
\item $f(x) = h(x)$ when $x$ is in the domain of $h$.
\end{itemize}
Adjacency is defined by restriction as 
in the proof of Proposition~\ref{prop:koenig}. 

The only difficulty of the proof is to show that
$\mathcal T$ has vertices on all levels $n$. 
We will first construct a homomorphism $p$
from $\Q.\Gamma[\{a_1,\dots,a_n\}]$ to $\Gamma$ with the property that $p(a_i) = a_i$
for $a_i$ in the domain of $h$, and
if $a_i - a_j = \infty$ for $i,j \leq n$,
then $p(a_i)$ and $p(a_j)$ are not $(h,q)$-connected.
Let $S$ be the set of points that are at distance at most $q$ of some $a_1,\dots,a_n$.
Let $S_1 \cup \cdots \cup S_k$ be the partition of $S$ induced by the copies of $\Z$ in $\Q.\Z$,
that is, $S_1,\dots,S_k$ are pairwise disjoint and each $S_i$ only contains points that lie in the same copy of $\Z$ in $\Q.\Z$.
Suppose without loss of generality that $S_1<\dots<S_{\ell-1}<S_\ell<S_{\ell+1}<\dots<S_k$ and that $S_\ell\subset\Z$, the standard copy in $\Q.\Z$.
For the elements $x \in S_\ell$ we set $p(x) := x$. 
For every $i\in\{1,\dots,k\}$, let $s_i$ and $t_i$ be the minimal and the maximal element of $S_i$, respectively.
Let $Q_\ell=\{z\in\Q.\Z\mid \exists z'\in S_\ell : |h(z')-z|\leq q\}$.
Write $S_\ell'$ for $h^{-1}(Q_\ell)$. 
If $\Z \setminus S_\ell' \cap \{z | z\geq t_i\}$ does not contain arbitrarily long intervals, then $\Gamma$ has a finite-range endomorphism by Lemma~\ref{lem:koenig2}, contrary to our assumptions.
So there exists a $z_\ell \in \Z$ such that $[z_\ell,z_\ell+t_{\ell+1}-s_{\ell+1}+2q] \cap S_\ell' = \emptyset$. 
For $x \in S_{\ell+1}$, we set $p(x) := x-s_{\ell+1}+z_\ell+q$.
As above, set $Q_{\ell+1}$ to be the set of points that are at distance at most $q$ of a point in $h(p(S_\ell\cup S_{\ell+1}))$.
Now, set $S_{\ell+1}':=h^{-1}(Q_{\ell+1})$. 
Then there exists a $z_{\ell+1} \in \Z$ such 
that $[z_{\ell+1},z_{\ell+1}+t_{\ell+2}-s_{\ell+2}+2q] \cap S_{\ell+1}' = \emptyset$. For $x \in S_{\ell+2}$, we set $p(x) := x-s_{\ell+2}+z_{\ell+1}+q$. 
Continuing in this way, we define $p$ for all $x \in \{a_1,\dots,a_m\}$ (the construction for $i<\ell$ is symmetric). The construction is illustrated in Figure~\ref{fig:lemma-6}.
We have that $p$ is a homomorphism $\Q.\Gamma[\{a_1,\dots,a_n\}]\to \Gamma$ since it is $\sim_q$-equivalent to the identity function of $\Q.\Z$.

\begin{center}
\begin{figure}
\begin{tikzpicture}
    \draw node at (0,0) {\includegraphics[width=1\linewidth]{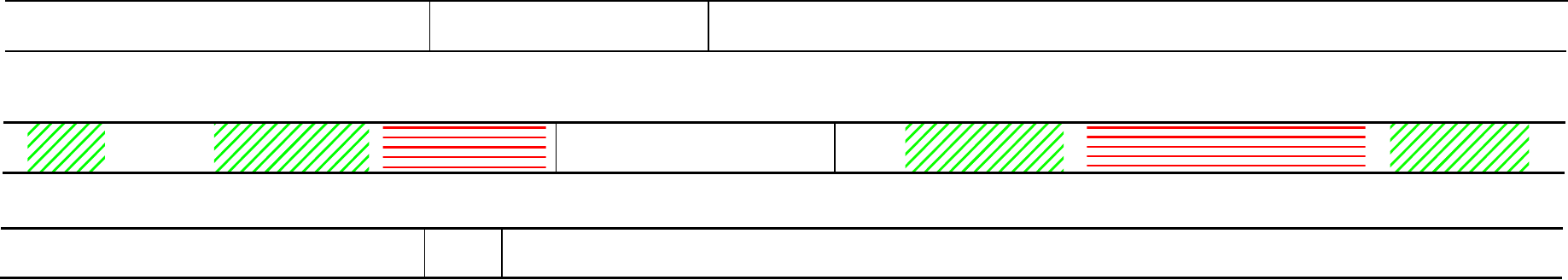}};
    \draw node (Sl) at (-0.9,-.1) {$S_\ell$};
    \draw node (Slplus) at (-1.9,1.05) {$S_{\ell+1}$};
    \draw node (Slmoins) at (-3.1,-1.55) {$S_{\ell-1}$};
    \draw node (pSlplus) at (4,.1) {};
    \draw node (pSlmoins) at (-6,-.3) {};

    \draw[->] plot[smooth] coordinates { ([yshift=-5]Slplus) (-0.9,0.5) (3,.5) ([yshift=3]pSlplus) };
    \draw[->] plot [smooth] coordinates { ([yshift=20]Slmoins) (-3.5,-.6) (-5.5,-.6) ([yshift=-2]pSlmoins) };
\end{tikzpicture}
\caption{Illustration of the construction of $p$ in Lemma~\ref{lem:hom-extension}. Each horizontal part is a copy of $\Z$, and the copies of $\Z$ are ordered from bottom to top. The part hashed diagonally is $h^{-1}(Q_\ell)$, i.e., the set of points that are mapped by $h$ to some point at distance $\leq q$ of $h(S_\ell)$. The part hashed horizontally is the set of points that are mapped by $h$
to some point at distance $\leq q$ from $h(p(S_{\ell+1}))$.}
\label{fig:lemma-6}
\end{figure}
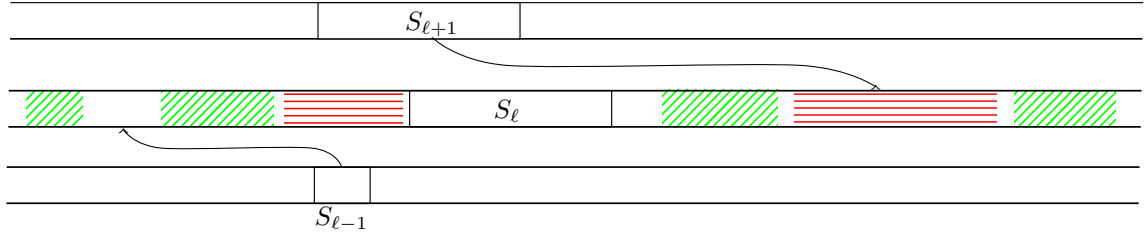
\end{center}


Let $r \colon \{a_1,\dots,a_n\} \to \Q.\Gamma$ be any map such that
\begin{itemize}
\item $r(x) = h(x)$ for all $x \in S_1$;
\item $r(a_i) - r(a_j) = (h\circ p)(a_i)-(h\circ p)(a_j)$ 
for all $i,j \leq n$, where either side of the equation can be $+\infty$ or $-\infty$.
\end{itemize}
Observe that by construction of $p$,
when $a_i - a_j = \infty$, then
$a_i, a_j$ are neither $(h \circ p,q)$-connected
nor $(r,q)$-connected. 
Hence, $r \sim_q (h \circ p)$, and
 Lemma~\ref{lem:substitution} implies that
 $r$ is a homomorphism from
$\Q.\Gamma[\{a_1,\dots,a_n\}]$
to $\Q.\Gamma$ which shows that
$\mathcal T$ has a vertex on level $n$. 

Observe that $\mathcal T$ is finitely branching, and by K\"onig's lemma contains an infinite branch $\mathcal B$.
Note that we cannot conclude with Proposition~\ref{prop:compactness}: we would only obtain a function which is $\sim_q$-equivalent to every
function in the infinite branch, but we also want a function that extends $h$.
However, a straightforward modification of the proof of Proposition~\ref{prop:compactness} using the hypothesis that the set $\{p\in\Q\mid \exists z\in\Z : h(z)\in p.\Q\}$ is bounded
yields the existence of the required function (in the second case of the inductive construction, $P$ and $Q$ are possibly infinite, but bounded, so that we can still choose a suitable rational $p$). 
\end{proof}

\begin{remark}
The hypothesis on the image of $h$ is necessary.
Indeed, let $\Gamma$ be simply $(\Z;<)$.
The function $h\colon\Z\to\Q.\Z$ which maps $z$ to $(z,0)$, i.e., that takes $z$ to the $0$ element of the $z$th copy of $\Z$ in $\Q.\Z$,
is a homomorphism $\Gamma\to\Q.\Gamma$.
However there is no endomorphism $\Q.\Gamma\to\Q.\Gamma$ that extends $h$. Such an endomorphism needs to map for example $(1,0)$
to an element in $\Q.\Z$ which is greater than every element in the image of $h$, but such an element does not exist.
\end{remark}


\subsection{Petrus}
\label{sect:petrus}
The following theorem is the rock upon which we build our church.

\begin{theorem}[Petrus ordinis]
\label{thm:petrus}
Let $\Gamma$ be a reduct of 
$({\mathbb Z};<)$ with finite relational signature and without an endomorphism of finite range. 
Then the following are equivalent:
\begin{enumerate}[(1)]
\item\label{item:equality} 
there exists a reduct $\Delta$ of $(\Q;<)$ such that $\Csp(\Delta)$ equals $\Csp(\Gamma)$;

\item\label{item:equality-endo} $\Q.\Gamma$ has an endomorphism whose range induces a structure isomorphic to a reduct of $(\Q;<)$;

\item\label{item:saturated-tight-rank} for all $t \geq 1$, there is an $e \in \End(\Q.\Gamma)$ and $z \in \Q.\Z$ such that 
$|e(z+t) - e(z)| > t$;

\item\label{item:unbounded-distance-degree} all binary relations with a primitive positive definition in $\Q.\Gamma$ are either empty, the equality relation, or have unbounded distance degree;

\item\label{item:split} for all distinct $z_1,z_2 \in \Z$ there is a homomorphism $h \colon \Gamma \to \Q.\Gamma$ such that $h(z_1) - h(z_2) = \infty$; 

\item\label{item:split-and-keep-separate} for all distinct $z_1,z_2 \in \Z$ there is an $e \in \End(\Q.\Gamma)$ such that $e(z_1) - e(z_2) = \infty$, and for
all $z_1,z_2 \in \Q.\Z$ with $z_1-z_2=\infty$ we have
$e(z_1)-e(z_2) = \infty$; 

\item\label{item:distribute} there exists an $e \in \End(\Q.\Gamma)$ with infinite range such that 
$e(x)-e(y) = \infty$ or $e(x)=e(y)$ for any two distinct $x,y \in \Q.\Gamma$. 
\end{enumerate}
\end{theorem}
\begin{proof}
Throughout the proof, let $q$ be the qe-degree of $\Gamma$, which is finite since $\Gamma$ has a finite signature. 

$(\ref{item:equality}) \Rightarrow (\ref{item:equality-endo})$. 
Since $\Delta$ has the same $\Csp$ as $\Gamma$, and $\Delta$
is $\omega$-categorical, a standard compactness argument (see e.g.~\cite{Bodirsky-HDR}) shows that
there is a homomorphism $f$ from the countable structure $\Q.\Gamma$ to $\Delta$. 
Lemma~\ref{lem:sat} asserts the existence of a homomorphism $g$ from $\Delta$ to 
$\Q.\Gamma$, 
because every pp-sentence
that is true in $\Delta$ is also true in $\Q.\Gamma$, and $\Q.\Gamma$ is saturated.
Then $e := g \circ f$ is an endomorphism of $\Q.\Gamma$. 

Since $\Delta$ is $\omega$-categorical, it has a model-complete core $\Delta'$~\cite{Cores-Journal}.
If $\Delta'$ is a finite structure, then $\Gamma$ has an endomorphism with a finite range,
contradicting our hypothesis. Therefore $\Delta'$ is an infinite structure which is again a reduct of $(\Q;<)$.
Let us henceforth assume that $\Delta$ is already a model-complete core.
$\Delta$ does not have non-injective endomorphisms since it is a core, and in particular $f\circ g$ is injective. 
Let $S=f(\Q.\Z)$. If $S$ were finite, then $f\circ g$ would be an endomorphism of $\Delta$
with finite range, in contradiction to the previous sentence. So $S$ is infinite.
$\Delta[S]$ is preserved by all the increasing functions, since all the increasing functions on $S$ can be extended to increasing functions on $\Q$,
which preserve $\Delta$. It follows that $\Delta[S]$ is isomorphic to a reduct of $(\Q;<)$.
Since $f\circ g$ is an endomorphism of $\Delta$, and since $\Delta$
is a core, $f\circ g$ is in fact an embedding. This implies that $g$ is an embedding,
so that $\Delta[S]$ is isomorphic to $\Q.\Gamma[g(S)]=\Q.\Gamma[(g\circ f)(\Q.\Z)]$, as required.

$(\ref{item:equality-endo}) \Rightarrow (\ref{item:saturated-tight-rank})$.
Let $f$ be an endomorphism of $\Q.\Gamma$ whose range induces in $\Q.\Gamma$ a structure isomorphic to a reduct of $(\Q;<)$.
Let $t\geq 1$. It is not possible that $f(z)=f(z+t)$ for all $z\in\Q.\Z$, for otherwise $\Gamma$ would have a finite-range endomorphism.
Indeed, we can restrict $f$ to a homomorphism $\Gamma\to\Q.\Gamma$ whose range is finite.
We can then construct a function $g\colon\Z\to\Q.\Z$ such that $f\sim_q g$ and such that the range of $g$ is contained in $\Z$.
This $g$ would then be an endomorphism of $\Gamma$ by Lemma~\ref{lem:substitution}, a contradiction.
Pick a $z\in\Q.\Z$ such that $f(z)\neq f(z+t)$.
The intersection of the image of $f$ with the set $\{z'\in\Q.\Z \mid |z'-f(z)|>t\}$ is nonempty, because this set is cofinite and we can assume that $f$ does not have a finite
range by the same argument as above. 
Pick a $z'\in\Q.\Z$ in the range of $f$ such that $|f(z)-z'|>t$; this is possible by the previous remark.
Since the range of $f$ is isomorphic to a reduct of $(\Q;<)$, we can in particular find an automorphism $\alpha$ of the structure induced by $\Q.\Gamma$ on the image of $f$
that maps $\{f(z),f(z+t)\}$ to $\{f(z),z'\}$. Therefore we have $|(\alpha\circ f)(z+t) - (\alpha\circ f)(z)| = |z'-f(z)|>t$. 

$(\ref{item:saturated-tight-rank}) \Rightarrow (\ref{item:unbounded-distance-degree})$. 
Let $R$ be a binary relation with a primitive positive definition in $\Q.\Gamma$.
Suppose that $R$ is not empty and is not the equality relation.
Let $k$ be the supremum of the integers $t$ such that there exists $(z_1,z_2)\in R$
with $|z_1-z_2|=t$. Since $R$ is neither empty nor the equality relation, it
follows that $k$ is positive. If $k$ is $\infty$, then $R$ has infinite distance degree.
Otherwise let $(z_1,z_2)$ be a pair in $R$ such that $|z_1-z_2|=k$.
Let $e$ be an endomorphism of $\Q.\Gamma$ and $z$ be such that $|e(z+k)-e(z)|>k$.
Let $\alpha$ be an automorphism of $\Q.\Gamma$ that maps $\{z_1,z_2\}$ to $\{z,z+k\}$.
Then $(e\circ\alpha)(z_1,z_2)$ is in $R$ since $R$ is preserved by the endomorphisms
of $\Q.\Gamma$ and by construction $|(e\circ\alpha)(z_1)-(e\circ\alpha)(z_2)|>k$, a contradiction to the choice of $k$.

$(\ref{item:unbounded-distance-degree}) \Rightarrow (\ref{item:split})$. 
Suppose that $(\ref{item:split})$ does not hold, that is,
there are distinct $z_1,z_2 \in \Z$ 
such that for all homomorphisms $h$ from $\Gamma$ to $\mathbb{Q}.\Gamma$ we have
that $h(z_1) - h(z_2) < \infty$. 
Then by Proposition~\ref{prop:koenig} there
is an $r \geq 0$ and a finite $S \subseteq \Z$ 
containing $\{a_1,a_2\}$ such that for all homomorphisms $f \colon \Gamma[S] \to \Gamma$ we have $|f(a_1) - f(a_2)| \leq r$. 
Now consider the following primitive positive formula $\phi$: the variables of $\phi$ are
the elements of $S$, all existentially quantified except $a_1$ and $a_2$, which are free. 
The formula $\phi$ contains the conjunct $R(x_1,\dots,x_n)$ for
a relation $R$ from $\Gamma$ if and only if
$\Gamma[S] \models R(x_1,\dots,x_n)$. 
Then $\phi$ defines a binary relation,
which has bounded distance degree by the previous discussion, and which is not the
equality relation since it contains the pair 
$(a_1,a_2)$. 

$(\ref{item:split}) \Rightarrow (\ref{item:split-and-keep-separate})$. 
An immediate consequence of Lemma~\ref{lem:hom-extension}. 

$(\ref{item:split-and-keep-separate}) \Rightarrow (\ref{item:distribute})$. 
Again an argument based on K\"onig's tree lemma. 
Let $a_1,a_2,\dots$ be an enumeration of $\mathbb{Q}.\Z$. 
Let $\mathcal T$ be a tree whose vertices
on the $i$-th level are the $\sim_q$-equivalence classes of homomorphisms $g$ from $\mathbb{Q}.\Gamma[\{a_1,\dots,a_n\}]$ to $\Gamma$ such that $g(a_i)-g(a_j) = \infty$ or $g(a_i) = g(a_j)$ for all $i,j \leq n$. Adjacency of vertices is defined by restriction between representatives. 
We have to show that the tree has vertices on all levels. 
Let $\{u_1,v_1\},\dots,\{u_k,v_k\}$ be an 
enumeration of all 2-element subsets of $\{a_1,\dots,a_n\}$. We will show by induction on $i \geq 0$ that there exists an endomorphism $f_i$ such that $f_j(u_j)-f(v_j) = \infty$ or $f_j(u_j) = f(v_j)$ for all $j \leq i$. 
The statement is trivial for $i=0$. 
So suppose we have already found $f_i$ for some $i \geq 0$, and want to find $f_{i+1}$. 
If $f_i(u_{i+1})-f_i(v_{i+1}) = \infty$ or $f_i(u_{i+1}) = f_i(v_{i+1})$
then there is nothing to show. Otherwise,
let $\alpha$ be an automorphism of $\mathbb{Q}.\Gamma$ that maps $f_i(u_{i+1})$ and $f_i(v_{i+1})$ to $\Z$. 
By~$(\ref{item:split-and-keep-separate})$, there exists an $e \in \End(\mathbb{Q}.\Gamma)$ such that $e(\alpha(f_i(u_{i+1})))-e(\alpha(f_i(v_{i+1}))) = \infty$,
and such that for all $x,y \in \Q.\Z$ we have that $e(x) - e(y) = \infty$. 
Hence, $f_{i+1} := e \circ \alpha \circ f_i$ has the desired property. The tree $\mathcal T$ has finitely many vertices on each level and hence must contain an infinite branch, which gives rise to an endomorphism of $\mathbb{Q}.\Gamma$ as in the proof of Proposition~\ref{prop:koenig}. 

$(\ref{item:distribute}) \Rightarrow (\ref{item:equality})$. Let $\Delta$ be the structure
induced by $\mathbb{Q}.\Gamma$ on the image of the endomorphism $e$
whose existence has been asserted in $(\ref{item:distribute})$. 
Note that a literal $x<y+k$ for $k\in\Z$ is true in $\Delta$ iff $x<y$ is true, and the same holds for literals $x\leq y+k$.
Therefore the restriction of the relations of $\Q.\Gamma$ to $\Delta$ are definable using $x<y$ and $x\leq y$.
Since $<$ is a dense linear order without endpoints on $\Delta$, it follows that $\Delta$ is isomorphic to a first-order reduct of $(\Q;<)$.
\end{proof}

\subsection{Boundedness and Rank}\label{sect:boundedness}
Let $\Gamma$ be a reduct of $(\Z;<)$ without a finite-range endomorphism. 
Theorem~\ref{thm:petrus} (Petrus)
characterized the ``degenerate case'' when $\Csp(\Gamma)$ is the CSP for a reduct of $(\Q;<)$.
For such $\Gamma$, as we have mentioned before, the complexity of the CSP has already been classified.
In the following we will therefore assume that the equivalent items of
Theorem~\ref{thm:petrus}, and in particular, item~(\ref{item:saturated-tight-rank}), do \emph{not} apply. To make the best use of those findings,
we introduce the
following terminology.

\begin{definition}
Let $k\in\N^+,c \in \N$.
A function $e \colon \kappa_1.\Z \to \kappa_2.\Z$ is \emph{$(k,c)$-bounded} if for all $u \in \kappa_1.\Z$ we have $|e(u+k) - e(u)| \leq c \, .$ 
\end{definition}
\vspace{-.5cm}
\begin{flushright}
\end{flushright}
\vspace{-.5cm}
We say that $e$ is \emph{tightly-$k$-bounded} if it is $(k,k)$-bounded,
and \emph{$k$-bounded} if it is $(k,c)$-bounded for some $c \in \N$. 
We say that $\kappa.\Gamma$ 
is \emph{$(k,c)$-bounded} if all its endomorphisms are; similarly, $\Gamma$ is \emph{tightly-$k$-bounded} if all its endomorphisms are. 
We call 
the smallest $t$ such that $\kappa.\Gamma$ is tightly-$t$-bounded the 
\emph{tight rank} of $\kappa.\Gamma$.
Similarly, we call the smallest $r$ such
that $\kappa.\Gamma$ is $r$-bounded 
the \emph{rank} of $\kappa.\Gamma$.
The negation of item~(\ref{item:saturated-tight-rank}) in Theorem~\ref{thm:petrus} says that
there exists a $t \in \mathbb N$ such that 
$\Q.\Gamma$ is tightly-$t$-bounded. 
Clearly, being tightly-$t$-bounded implies
being $t$-bounded. Hence,
the negation of item~$(\ref{item:saturated-tight-rank})$ in Theorem~\ref{thm:petrus} also implies
that $\Q.\Gamma$ has finite rank $r \leq t$.

\begin{example}\label{ex:rank1}
There are rank one reducts of $(\Z;<)$
which do have non-injective endomorphisms,
but no finite-range endomorphisms. Consider the second structure in the Example~\ref{ex:F}:
$$\Gamma := (\Z; \Diff{\{2\}}, \{(x,y) : |x-y| \leq 2\}) \; .$$
Note that 
$\Gamma$ has rank one: as every endomorphism $e$ preserves the
relation $\{(x,y) : |x-y| \leq 2\})$ we have
$|e(x+1) - e(x)| \leq 2$. 
Also note that $\Gamma$ 
has the non-injective endomorphism $e$ defined by 
$e(x) = x$ for even $x$, and $e(x) = x+1$ for odd $x$. 
\end{example}





Sections~\ref{sect:rank-one} and~\ref{sect:rank-greater-one} are devoted to proving that one
can replace $\Gamma$ by another reduct $\Delta$ of $(\Z;<)$ which has the same CSP
and such that $\Q.\Delta$ has both rank one and tight-rank one.
We first prove a general fact that holds in both cases.

\begin{lemma}\label{lem:rank-one}
Let $\Gamma$ be a reduct of $(\Z;<)$ with finite relational signature
such that $\Q.\Gamma$ has rank $r$ (resp. tight-rank $t$). Then $\Gamma$ has rank $r'$ with $r'\leq r$ (resp. tight-rank $t'\leq t$). 
\end{lemma}

\begin{proof}
Suppose that for all $c \in \N$ there is
$f \in \End(\Gamma)$ such that $|f(z+r)-f(z)| > c$ for some $z \in \Z$. Viewing
$f$ as a homomorphism from $\Gamma$ to $\Q.\Gamma$, we can apply Lemma~\ref{lem:hom-extension} and
find an endomorphism $e \in \End(\Q.\Gamma)$ which extends $f$; this endomorphism shows that $\Q.\Gamma$ is not $(r,c)$-bounded. 
Similarly, an endomorphism $f$ of $\Gamma$ so that $|f(z+t)-f(z)|>t$ for some $z\in\Z$
extends to an endomorphism of $\Q.\Gamma$.
\end{proof}

\subsubsection{The Rank One Case}
\label{sect:rank-one}

The main result of this section, Theorem~\ref{thm:rank-one}, implies that for each rank one reduct $\Gamma$ of $(\Z;<)$ without finite range endomorphisms
 there exists a reduct $\Delta$ of $(\Z;<)$ which has the same CSP as $\Gamma$ and where $\suc$ is pp-definable, or for all $k \geq 1$ the relation $\Dist{k}$ is pp-definable.
Our strategy is the following: we study some properties of the endomorphisms of the rank one reduct $\Gamma$ of $(\Z;<)$, in the view of using Theorem~\ref{thm:galois} for $\End(\Q.\Gamma)$.
In general, the endomorphisms of $\Gamma$ are better behaved than the endomorphisms of $\Q.\Gamma$, as the latter endomorphisms
can exhibit different behaviours in each copy of $\Z$, and can collapse copies, whereas the former endomorphisms are more uniform, as we will show below.
Theorem~\ref{thm:rank-one-factor} is the first milestone in the strategy, as it allows us to replace $\Gamma$ with a reduct $\Delta$ of $(\Z;<)$
such that the endomorphisms of $\Q.\Delta$ are well-behaved, and thus are easier to understand.
The final step towards proving Theorem~\ref{thm:rank-one} is the classification of automorphisms groups of rank one reducts of $(\Q.\Z;<)$ in Theorem~\ref{thm:autos}.

\begin{lemma}\label{lem:little}
Let $e \colon \Z \to\Z$
be tightly-$t$-bounded and $(1,c)$-bounded for some $c,t\in \N$.
Then for all $n \in \N$, and $z \in \Z$, $|e(z+n) - e(z)| \leq n + ct \, .$
\end{lemma}
\begin{proof}
Let $n = pt + q$ for $0 \leq q \leq t$.
We have $|e(z+pt+q) - e(z+pt) | \leq qc$ by $q$ applications of $(1,c)$-boundedness,
and $|e(z+pt) -e(z)| \leq p t$ by $p$ applications of tight rank $t$.
We obtain
\begin{align*} 
& |e(z+n) - e(z)| \\
\leq \; & |e(z+pt+q) - e(z+pt)| + |e(z+pt) - e(z)|  \\
 \leq \; & qc + pt \\
 = \; & n+c(q-1)  \\
 \leq \; & n+ct 
\end{align*}
by the triangle inequality.
\end{proof}

The following can be shown by the same proof 
as the proof of Lemma 6 in~\cite{BodDalMarMotPin};
since our statement is more general, and since 
we use 
rank and tight rank instead of bounded distance degree, we still give the proof here for the convenience of the reader. We note that this lemma must be proved in the atomic model.

\begin{lemma}
\label{lem:6}
Let $e\colon\Z\to\Z$ be tightly-$t$-bounded and $(1,c)$-bounded.
Then either $e$ together with translations locally generates a
function with finite range, or there exists $k > ct+1$ such that for all $x,y \in \Z$ with $|x-y|=k$ we have $|e(x)-e(y)| \geq k$. 
\end{lemma}
\begin{proof}
Assume for all $k > ct+1$ there are $x,y \in \Z$ with $|x-y|=k$ and $|e(x)-e(y)| < k$. 
We will prove that $e$ locally generates a function with range of size at most $2ct+1$.
By a by now standard application of K\"{o}nig's tree lemma,
it suffices to show that for every finite $A\subseteq \Z$,
$e$ generates an $h$ such that
$|h(A)| \leq 2ct +1$. 

    Enumerate the pairs $(x,y) \in A^2$ with $x<y$ by $(x_1,y_1),\ldots,(x_r,y_r)$. 
       Let $m$ be the smallest number with the property that 
       $\mathcal F := \{e\} \cup \Aut(\Z;<)$ generates an endomorphism $h_1$ 
       such that $|h_1(x_1)-h_1(y_1)| = m$. 
    We claim that $m \leq ct + 1$. 
	Otherwise, by assumption there are $x,y \in \Z$
	with $|x-y|=m$ and $|e(x)-e(y)| < m$.
	Let $a$ be the automorphism of $(\Z;<)$
	such that $a(\{h_1(x),h_1(y)\}) = \{x_1,y_1\}$. 
	Then $\mathcal F$ also generates $h_1' := e \circ a \circ h_1$,
	but $|h_1'(x_1) - h_1'(y_1)| < m$ in contradiction to the choice of $m$.  We conclude that $\Gamma$
	has an endomorphism $h_1$ such that $|h_1(x_1)-h_1(y_1)| \leq ct +1$.
    
    Similarly, there exists $h_{2}$ generated by $\mathcal F$ such that $|h_{2}(h_1 (x_{2})) - h_{2}(h_1 (y_{2}))| \leq ct +1$. Continuing like this we arrive at a function $h_{r}$ generated by $\mathcal F$ such that $$|h_{r}h_{r-1}\cdots h_1 (x_{r})-h_{r}h_{r-1}\cdots h_1 (y_{r})|\leq ct +1.$$ Now consider $h := h_r\circ\cdots \circ h_1$. 
    Set $f_j:=h_{r}\circ \cdots\circ h_{j+1}$ and $g_j:=h_{j}\circ\cdots\circ h_{1}$, for all $1\leq j\leq r$; so $h=f_j\circ g_j$. Then, since by construction $|g_j(x_j)-g_j(y_j)|\leq ct +1$, we have that  for all $j \in \Z$ with $1\leq j \leq r$,
    \begin{align*}
    & |h(x_j)-h(y_j)| \\
    = & \; |f_j(g_j(x_j))-f_j(g_j(y_j))| \\
    \leq & \; |g_j(x_j)-g_j(y_j)|+ct && \text{(Lemma~\ref{lem:little})} \\
    \leq & \; 2 c t+1 \, ,
    \end{align*}
    and our claim follows.
\end{proof}


\begin{definition}\label{def:stable}
Given $e \colon \kappa_1.\Z \to \kappa_2.\Z$, we call $s \in \N^+$ 
\emph{stable for $e$} if $|e(z+s)-e(z)|=s$ for all $z\in\kappa_1.\Z$.
\end{definition}

\begin{remark}
	In the locally finite case~\cite{BodDalMarMotPin}, a number $s$ was defined to be stable for $e\colon\Z\rightarrow\Z$
	if for all $z\in\Z$, $e(z+s)=e(z)+s$ or for all $z\in\Z, e(z+s)=e(z)-s$.
	The definition used in this paper is strictly weaker, but the previous definition
	is too strict for functions $\Q.\Z\rightarrow\Q.\Z$.
	For example, the function $e$ that maps $x$ to $-x$ on one copy of $\Z$,
	and that is the identity on the other copies has no stable number if we consider
	the natural generalization of the old definition, whereas for our purposes it is easier
	to consider that $1$ is stable for this function.
\end{remark}

\begin{lemma}
\label{lem:tight-rank-stable}
Let $e\colon\Z\to\Z$ be tightly-$t$-bounded and $1$-bounded.
Then $t$ is stable for $e$, or $e$ locally generates with translations a function with finite range.
\end{lemma}
\begin{proof}
Let $c\in\N$ be such that $e$ is $(1,c)$-bounded, and assume that $e$ does not locally generate a function with finite range.
By Lemma~\ref{lem:6}, there exists $k > ct+1$ so that for all $z$ we have
$|e(z+k)-e(z)| \geq k$,
and hence either $e(z+k) \geq e(z)+k$ 
or $e(z+k) \leq e(z)-k$ 
for each $z\in \Z$. 
We will first show that either 
$e(z+k) \geq e(z)+k$ 
for all $z\in \Z$,
or $e(z+k) \leq e(k)-k$ 
for all $z\in \Z$. 
Suppose otherwise that 
there are $z_1,z_2 \in \Z$
such that $e(z_1+k) \geq e(z_1)+k$
and $e(z_2+k) \leq e(z_2) - k$.
Clearly, we can choose $z_1,z_2$ such that $|z_1-z_2|=1$. We only treat the case
that $z_2 = z_1 + 1$, since the other case is symmetric. Then
\begin{align*}
e(z_2) - e(z_2+k) & \geq k && \text{by assumption,} \\
- e(z_2) + e(z_1) & \geq - c  && \text{by 1-boundedness,} \\
e(z_2+k) - e(z_1+k) & \geq - c && \text{by 1-boundedness,} \\
e(z_1+k) - e(z_1) & \geq k  && \text{by assumption.}
\end{align*}
Summing over those inequalities yields $ 0 \geq 2k - 2c \; ,$
a contradiction since $k>c$. 

In the following we assume without loss of generality that $e(z+k) \geq e(z)+k$ 
for all $z\in \Z$. 
Recall that 
$|e(z+t)-e(z)| \leq t$ for all $z \in \Z$ because $e$ is tightly-$t$-bounded. 
We next claim that $e(z+kt)= e(z)+kt$ for all $z \in \Z$. Since 
points at distance $t$ cannot be mapped to points at larger distance, we get that $e(z+kt) - e(z) \leq kt$. On the other hand, since 
$e(z+k) \geq e(z)+k$ for all $z \in \Z$, we obtain that $e(z+kt) \geq e(z)+kt$, proving the claim.

We now show that $|e(z+t) - e(z)| \geq t$ for all $z \in \Z$. Note that
\begin{align*} e(z) + k t \; = \; & e(z+k t) \\
= \; & e(z+t+ (k-1)t) \\
\leq \; & e(z+t)+(k-1)t
\end{align*} 
the latter inequality holding since $e(z+mt)-e(z) \leq m t$ for each $m \in \mathbb N$. Subtracting $(k-1)t + e(z)$ on both sides, our claim follows.
Since $|e(z+t)-e(z)| \leq t$ for all $z\in\mathbb{Z}$,
we obtain that $e(z+t) - e(z) = t$ and 
have proved the lemma.
\end{proof}

\begin{corollary}\label{cor:existence-stable-saturated}
 Let $\Gamma$ be a reduct of $(\Z;<)$ without finite range endomorphism
 so that $\Q.\Gamma$ has tight-rank $t$ and rank one.
 Then some integer is stable for every endomorphism $e\in\End(\Q.\Gamma)$.
\end{corollary}
\begin{proof}
 By Lemma~\ref{lem:rank-one}, $\Gamma$ has tight rank $t'\leq t$ and rank one.
 Let $z\in\Q.\Z$ and $e\in\End(\Q.\Gamma)$.
 Since $\Q.\Gamma$ has rank one, we have $e(z+k)-e(z)<\infty$ for all $k\in\Z$.
 As a consequence, $e$ induces an endomorphism $e'\colon\Gamma\to\Gamma$ by restricting $e$
 to the copy of $\Gamma$ that contains $z$.
 By Lemma~\ref{lem:tight-rank-stable}, $t'$ is stable for $e'$, so that $|e'(z+t')-e'(z)|=t'$.
 Since $e'(z)=e(z)$ and $e'(z+t')=e(z+t')$, $t'$ is stable for $e$.
\end{proof}

\begin{lemma}
\label{lem:stable-multiples}
    Let $e$ be a function from $\Z$ to $\Z$ and suppose that $s$ is stable for $e$ such that $s$ is minimal. Then the stable numbers for $e$ are precisely the multiples of $s$. \end{lemma}
\begin{proof}
Clearly, all multiples of $s$ are stable for $e$. 
For the other direction suppose that $p$
is stable but not divisible by $s$.

Write $p = ms + r$ where $m,r$ are positive integers and $0 < r < s$. Since $r$ is not stable
there exists $z \in \Z$ such that 
$e(z+r)-e(z) \neq r$.
But this is impossible since 
\begin{align*}
e(z+r)-e(z) = & \; e(z+p-ms)-e(z) \\
= & \; e(z-ms)+p-e(z) \\
= & \; e(z)-ms+p-e(z) \\
= & \; r\;.
\end{align*}
\end{proof}

\begin{lemma}\label{lem:infiniteCore}
    Let $e\colon\Z\to\Z$ be tightly-$t$-bounded and $1$-bounded,
    and let $s$ be stable for $e$ such that $s$ is minimal. 
    Then $\{e\} \cup \Aut(\Z;<)$ generates an $f$ 
    such that $f(\Z)=\{s\cdot z : z\in{\mathbb Z}\}$.
\end{lemma}
\begin{proof}
Assume that $s>1$ (if $s=1$ then $f$ can be chosen to be the identity and there is nothing to do).
Let $M$ be the monoid generated by $\{e\}\cup\Aut(\Z;<)$.
We claim that there exists an $f_1 \in M$ such that 
$f_1(0)=0$ and 
$f_1(1) \in \{s\cdot z : z\in \Z \}$. To see this, observe that since $s>1$ 
there exist $x_0,y_0\in \Z$ with $|x_0-y_0|=1$ 
and $|e(x_0)-e(y_0)|>1$. 
Write $r_1:=|e(x_0)-e(y_0)|$. If $r_1$ is not a multiple of $s$, then by Lemma~\ref{lem:stable-multiples} and Lemma~\ref{lem:tight-rank-stable}, $e$ is not tightly-$r_1$-bounded
and 
there exist $x_1,y_1\in\Z$ with $|x_1-y_1|=r_1$ and $|e(x_1)-e(y_1)|=:r_2>r_1$. 
Again, if $r_2$ is not a multiple of $s$, then there exist $x_2,y_2\in\Z$ with $|x_2-y_2|=r_2$ and $|e(x_2)-e(y_2)|=:r_3>r_2$. Continuing in this way, we obtain a sequence $(x_i,y_i)$ of pairs of distance $r_i$ (setting $r_0:=1$). By exchanging $x_{i+1}$ and $y_{i+1}$ if necessary, we may assume that $x_{i+1}<y_{i+1}$ iff $e(x_i)< e(y_i)$ for all $i$. There exist automorphisms $\alpha_i$ of $(\Z;<)$ such that $\alpha_i(e(x_i)) = x_{i+1}$ and $\alpha_i(e(y_i))=y_{i+1}$. Set $g_i := \alpha_i \circ e \circ \alpha_{i-1}\circ\cdots\circ \alpha_0\circ e$. Then $g_i \in M$, and $g_i$ sends $(x_0,y_0)$ to $(x_{i+1},y_{i+1})$, a pair of distance $r_{i+1}>r_i>\cdots>r_0$. Recall
$|e(x_0)-e(y_0)| \leq c$ since $e$ is $(1,c)$-bounded, so the sequence must end at some finite $i$. By construction, this happens only if $r_{i+1}$ is a multiple of $s$. Therefore, $r_{i+1}=|g_i(x_0)-g_i(y_0)| \in \{s\cdot z : z\in\Z\}$. 
By composing with an automorphism of $(\Z;<)$ we may assume that $x_0=0$, $y_0=1$, and $g_i(0)=0$. Set $f_1:=g_i$.

If $s>2$, then consider the number $f_1(2)$. We claim that there is $f_2 \in M$ such that $f_2(0)=0$ and $f_2(f_1(2))$ is a multiple of $s$. If already $f_1(2)$ is a multiple of $s$, then we can choose $f_2$ to be the identity. Otherwise, we can increase the distance between $f_1(2)$ and $0$ successively by applying shifts and $e$ just as before, where we moved away $1$ from $0$. After a finite number of steps, we arrive at a function $f_2$ such that $|f_2(0)-f_2(f_1(2))|$ is a multiple of $s$. Applying a shift one more time, we may assume that $f_2(0)=0$, and so $f_2$ has the desired properties.

We continue inductively, constructing for every $i<s$ a function $f_i$ such that $f_i(0)=0$ and $f_i\circ \cdots \circ f_1(i)$ is a multiple of $s$.
At the end, we set $f:=f_{s-1}\circ\cdots\circ f_1$. Since $s$ is stable for $e$, it is also stable for
 $f$, as $f$ is composed from $e$ and automorphisms of $(\Z; <)$. 
Hence, $f(\Z)$ contains the set $\{s \cdot z : z\in\Z\}$. For the other inclusion, let $v\in \Z$ be arbitrary, and write
$v= s \cdot z+ r$, where $z\in \Z$ and $0\leq r<s$. Then $f(s\cdot z+r)-f(r)$ is a multiple of $s$ since $s$ is stable for $f$.
By construction, $f(r)$ is a multiple of $s$ as well, so that $f(v)\in\{s\cdot z \mid z\in \Z\}$.
\end{proof}

The following definition arises naturally from the statement of Lemma~\ref{lem:infiniteCore}. 
\begin{definition}
Let $\Gamma$ be a structure over $\Z$ and let $k \in \N^+$. Then we write $\Gamma / k$
for the substructure of $\Gamma$ induced by the set $\{z \in \Z : z=0 \bmod k\}$. 
\end{definition}

\begin{lemma}\label{lem:division}
For all reducts $\Gamma$ of $(\Z;<)$ and $k\in\N^+$, the structure $\Gamma / k$
is isomorphic to a reduct of $(\Z;<)$, the isomorphism being the function $x\mapsto x/k$.
\end{lemma}
\begin{proof}
Let $R$ be an $n$-ary relation of $\Gamma$, and let $\phi$ be a standardised formula defining $R$ over $(\Z;<)$. 
Construct a formula $\phi'$ as follows: For all $i \in \Z$, replace every atomic formula of the form $x \leq y+i$ by $x\leq y+ \lfloor i/k\rfloor$,
and similarly for formulas of the form $x<y+i$.
We prove by structural induction on $\phi$ that for all
$z_1,\dots,z_n \in \Gamma/k$ we have $(\Z,<) \models \phi(z_1,\dots,z_n) \,\Leftrightarrow \, (\Z;<) \models \phi'(z_1/k,\dots,z_n/k)$.
If $\phi$ is $x \leq y+i$ for some $i \in \Z$, then $\Gamma/k\models\phi(x,y)$ iff $x\leq y+i$ iff $x/k \leq y/k + \lfloor i/k\rfloor$.
The case where $\phi$ is $x<y+i$ is treated similarly.
The cases of conjunction, disjunction, and negation follow immediately from the induction hypothesis.\end{proof}

For instance, in Example~\ref{ex:rank1} 
the structure $\Gamma/2$ is isomorphic to $(\Z;\suc,\{(x,y) : |x-y| \leq 1\})$.
If $\Delta$ is the reduct of $(\Z; <)$ isomorphic to $\Gamma/k$, we have that $\Q.\Delta$ is also
isomorphic to $\Q.(\Gamma/k)$, an isomorphism being $(q,z)\mapsto (q,k\cdot z)$, where $(q,z)\in\Q.\Z$.

\begin{theorem}\label{thm:rank-one-factor}
Let $\Gamma$ be a finite-signature reduct of 
$(\Z;<)$ without finite range endomorphisms and such that $\Q.\Gamma$ has rank one. 
Then $\Gamma$ has an endomorphism that maps $\Gamma$ to $\Gamma/k$ for some $k \in \N^+$,
and so that $\Q.(\Gamma/k)$ has tight-rank one.
\end{theorem}
\begin{proof}
Let $t$ be the tight-rank of $\Q.\Gamma$,
and let $c$ be such that $\Q.\Gamma$ is $(1,c)$-bounded.
Note in particular that every homomorphism $e\colon\Gamma\rightarrow\Q.\Gamma$ is
also $(1,c)$-bounded, since by Lemma~\ref{lem:hom-extension} it can
be extended to an endomorphism of $\Q.\Gamma$.
By Lemma~\ref{lem:rank-one}, $\Gamma$ has tight-rank $t'$, with $t'\leq t$.
By Corollary~\ref{cor:existence-stable-saturated}, every endomorphism of $\Q.\Gamma$
has a stable number, and in particular each endomorphism has a minimal one.
If the minimal stable number of every endomorphism is $1$, then $\Q.\Gamma$ has tight rank one
and we are done, choosing $k=1$.
Otherwise there exists an $e \in \End(\Q.\Gamma)$ 
such that $1$ is not stable.
Since $\Q.\Gamma$ has rank one, $e$ sends copies of $\Z$ of copies of $\Z$.
Moreover $1$ is not stable for $e$, so there exists a copy of $\Z$ and some integer $s>1$ such that
$s$ is stable for the restriction of $e$ to that copy, which we call $\hat e$, and so that no $s'$ with $s'<s$
is stable for $\hat e$.
By Lemma~\ref{lem:infiniteCore}, there
exists a function $f$ generated by $\{\hat e\} \cup \Aut(\Z;<)$ such
that $f(\Z) = \{s \cdot z : z \in \Z\}$.
By Lemma~\ref{lem:tight-rank-stable}, $t'$ is stable for $f$,
and $t'$ is divisible by $s$ since $|f(z+t')-f(z)|=t'$ and $f(z+t'),f(z)\in\{s\cdot z : z\in\Z\}$.

Observe that $\Gamma/s$ cannot have a finite range endomorphism: if $g$ were such an endomorphism, then $g \circ f$ would be a finite range endomorphism for $\Gamma$, contrary to our assumption. 
By Lemma~\ref{lem:division}, $\Gamma/{s}$ is isomorphic to a reduct $\Delta$ of $(\Z;<)$.
It is also clear that the function $(a,z)\mapsto (a, sz)$ from $\Q.\Z$ to $\Q.\Z$
is a homomorphism between $\Q.\Delta$ and $\Q.\Gamma$.
We claim that $\Q.\Delta$ has rank one and tight-rank at most $t'/s$.

Let $e\in\End(\Q.\Delta)$.
Let $x\in\Q.\Z$, and $h$ be the embedding $\Delta\rightarrow\Q.\Delta$
that maps $\Z$ to the copy of $\Z$ that contains $x$ in $\Q.\Z$.
Note that $1$ is stable for $h$.
Define $e'(z) = s \cdot e(h(f(z)/s))$, which is a homomorphism $\Gamma\rightarrow\Q.\Gamma$.
As such, $e'$ is $(1,c)$-bounded. 
Since $f$ is surjective as a function $\Z\rightarrow\{sz\mid z\in\Z\}$, there exists $y\in\Z$ such that $h(f(y)/s) = x$.
Then
\begin{align*}
	|e(x+1) - e(x)| & =  \left|e\left(h\left(\frac{f(y)}{s}\right) + 1\right) - e\left(h\left(\frac{f(y)}{s}\right)\right)\right| \\
				&  = \left|e\left(h\left(\frac{f(y)+s}{s}\right)\right) - e\left(h\left(\frac{f(y)}{s}\right)\right)\right|\\
			      &= \frac{1}{s}\cdot |e'(y\pm s) - e'(y)| \leq c,
\end{align*}
where the last inequality holds by $(1,c)$-boundedness of $e'$ applied $s$ times.
Hence, all the endomorphisms of $\Q.\Delta$ are $(1,c)$-bounded and $\Q.\Delta$ has rank one. Similarly, we have
\begin{align*}
	\left|e\left(x+\frac{t'}{s}\right) - e(x)\right| &= \left|e\left(h\left(\frac{f(y)}{s}\right) + \frac{t'}{s}\right) - e\left(h\left(\frac{f(y)}{s}\right)\right)\right| \\
			      &= \frac{1}{s}\cdot |e'(y\pm t') - e'(y)| \leq \frac{t'}{s},
\end{align*}
i.e., $e$ is tightly-$t'/s$-bounded and $\Q.\Delta$ has tight rank at most $t'/s$.

Since $\Delta$ satisfies all assumptions that we had on $\Gamma$, we may repeat the argument.
If all endomorphisms of $\Delta$ are automorphisms, then we are done. 
This process terminates, since the tight rank of $\Q.\Delta$ is bounded above by $t'/s$, which
is strictly smaller than the tight rank of $\Q.\Gamma$.
Observe furthermore that if $\Delta'$ is the reduct of $(\Z;<)$ that is isomorphic to $\Delta/{s'}$, then $\Delta'$ is isomorphic to $\Gamma/{ss'}$ by the obvious composition of isomorphisms,
so that the resulting structure at termination is indeed of the form $\Gamma/k$ for some $k \in \N$. 
\end{proof}

\begin{lemma}\label{lem:injective-surjective}
    Let $\Gamma$ be a reduct of $(\Z;<)$ with rank one and finite tight-rank, whose automorphisms contain the translations,
    and suppose that $\Gamma$ does not have finite range endomorphisms. 
    Then every injective or surjective endomorphism of $\Gamma$ is an isometry and an automorphism of $\Gamma$.  
\end{lemma}
\begin{proof}
Let $e \in \End(\Gamma)$ be injective (or surjective). 
By Lemma~\ref{lem:tight-rank-stable}, there exists
a stable number $s$ for $e$; choose $s$ such that
$s$ is minimal. 
By Lemma~\ref{lem:infiniteCore}, $\{e\} \cup \Aut(\Z;<)$ generates a function $f$ such that
$f(\Z) = \{s z : z \in \Z\}$. 
Suppose that $k$ is stable for $f$, and choose $k$
so that $k$ is minimal. 
We claim that $k = s$. 
Since $|f(x+k) - f(x)|$ is a multiple of $s$, we
have $|f(x+k) - f(x)| \geq s$ and hence $k \geq s$.
Moreover, 
$|f(x+s)-f(x)| = s$
since $f$ is a composition of $e$ and automorphisms of $(\Z; <)$. 
So $s$ is stable for $f$, and so $k \leq s$. 

Since $e$ is injective (or surjective), the function $f$
is also injective (or surjective, respectively). 
This is only possible when $s=1$. 
To see this, suppose without loss of generality that
$f(0)=0$ (otherwise, apply the argument below to the injective endomorphism $x \mapsto x - f(0)$ instead of $f$).
We have $f(1) = z_0 s$ for some $z_0 \in \Z \setminus \{0\}$
since $f(\Z) = \{sz : z\in\Z\}$.
Moreover, $f(f(1)) = f(z_0 s) \in \{z_0 s,-z_0 s\}$ since $s$ is stable for $f$. 
If $f(z_0 s) = z_0 s$ and when $f$ is injective,
 then this implies that $z_0 s = 1$ and hence
 that $|s| = 1$. Otherwise, $f(z_0 s) = -z_0s$. 
Observe that $f(-z_0 s) \in \{z_0 s,-z_0 s\}$. 
If $f(-z_0s) = -z_0s$ then $f$ cannot
be injective since $-z_0s$ and $z_0s$ are distinct but have equal function values. If $f(-z_0s) = z_0s$
then the injectivity of $f$ implies that $-z_0s = 1$,
and hence that $s=1$. 

When $f$ is surjective, then $\{z  s : z \in \Z \} = \Z$, and $s = 1$.

So we conclude that $s=1$, 
and that $e$ is either an automorphism of
$(\Z;<)$, or the composition of $x \mapsto -x$ with such an automorphism, and hence isometric. 
In both cases, there exists an $\alpha \in \Aut(\Z;<)$ such that $\alpha \circ e \circ e$ is the identity, 
so $e$ is invertible in $\End(\Gamma)$ and hence
an automorphism of $\Gamma$. 
\end{proof}


Let $i \colon \Q.\Z \to \Q.\Z$
be defined by $x \mapsto -x$. 
Note that $i$ fixes the copies 
of $\Z$ in $\Q.\Z$ setwise. 

\begin{theorem}
\label{thm:autos}
Let $\Gamma$ be a reduct of 
$(\Z;<)$ without finite range endomorphisms and such that $\Q.\Gamma$ has rank one.
Then the automorphism group of $\Q.\Gamma$ 
equals exactly one of the following: 
\begin{itemize}
\item $\Aut(\Q.\Z; \Dist1)$;
\item $\Aut(\Q.\Z;F)$  (Example~\ref{ex:F}), which is generated by $\Aut(\Q.\Z;<)$ and $i$.
\item $\Aut(\Q.\Z;<)$.
\end{itemize}
\end{theorem}
\begin{proof}
Since $\Q.\Gamma$ is a reduct of $(\Q.\Z;<)$ we have $\Aut(\Q.\Z;<) \subseteq \Aut(\Q.\Gamma)$. 
Since $i$ does not preserve $<$, 
we have that $\Aut(\Q.\Z;<) \subsetneq \Aut(\Q.\Z;F)$. 

We claim that $\Aut(\Q.\Gamma)$ is
contained in $\Aut(\Q.\Z; \Dist1)$. 
Let $\alpha$ be an arbitrary automorphism
of $\Q.\Gamma$. Since $\Gamma$ has rank one, we have that $x-y=\infty$ if and only if
$\alpha(x)-\alpha(y) = \infty$. In other words,
$\alpha$ permutes the copies of $\Z$ in $\Q.\Z$.
Moreover, when we identify $\Z$ and $\alpha(\Z)$,
then the restriction of $\alpha$ to $\Z$ is
an automorphism of $\Gamma$; this holds for
all copies of $\Z$ in $\Q.\Z$.
As a direct consequence of Lemma~\ref{lem:injective-surjective} we have 
$\Aut(\Gamma) \subseteq \Aut(\Z;\Dist1)$. Hence, $\alpha$ is an automorphism of $(\Q.\Z;\Dist1)$, and this proves the claim. 

Clearly, $i$ preserves $\Dist1$, 
and therefore $\Aut(\Q.\Z;F) \subseteq \Aut(\Q.\Z; \Dist1)$. On the other hand,
$\Aut(\Q.\Z;F)$ also preserves the relation
$F$, and this relation is not preserved
by the automorphism of $(\Q.\Z;\Dist1)$ which acts as $x \mapsto -x$ on only one copy, and acts as the identity on all other copies.
Hence, $\Aut(\Q.\Z;F) \subsetneq \Aut(\Q.\Z; \Dist1)$. 

Next, we prove that if $\Aut(\Q.\Gamma)$ contains a permutation $\alpha$ not from $\Aut(\Q.\Z;<)$ 
then $\Aut(\Q.\Gamma)$ contains $i$.
Let $q_1,q_2,\dots$ be an enumeration of $\Q$.
Since $\alpha$ is not from $\Aut(\Q.\Z;<)$,
there exists a copy of $\Z$ in $\Q.\Z$ such that the restriction of $\alpha$ to
this copy does not preserve $<$; without loss of generality, 
we identify this copy with $q_1.\Z$. 
Since $\Q.\Gamma$ has rank one,
by composing $\alpha$ with an automorphism
of $(\Q.\Z;<)$, we can assume without loss of generality that
$\alpha$ fixes all copies of $(\Z;<)$ in $(\Q.\Z;<)$ setwise,
and $\alpha(0)=0$ in each of those copies. 
Since $\alpha$ preserves $\Dist1$, the restriction of $\alpha$ to $q_1.\Z$ equals $x \mapsto -x$. 
If the restriction of $\alpha$ to infinitely many other copies also has this form, then $i$ clearly lies in the closure of $\{\alpha\} \cup \Aut(\Q.\Z;<)$.
Otherwise, $\alpha$ equals the identity on all but finitely many copies. Without loss of generality, we can assume that $M$ is such that 
$\alpha$ acts as $i$ on the copies $q_1.\Z,\dots,q_M.\Z$ and is the identity on the copies $q_n.\Z$ with $n>M$.
Let $\beta_n$ be an automorphism of $(\Q.\Z;<)$ that maps $q_n.\Z$ to $q_1.\Z$ and for each $m<n$,
$\beta_n$ maps $q_m.\Z$ to some $q_{m'}.\Z$ with $m'>n$. 
Then $[\beta_k^{-1},\alpha]\circ\dots\circ[\beta^{-1}_{M+1},\alpha]$ (with $[f,g]:= f\circ g\circ f^{-1}\circ g^{-1}$ denoting the commutator)
acts as $i$ on $q_1.\Z \cup \cdots \cup q_k.\Z$,
and hence $i$ lies in the closure of the maps generated by $\{\alpha\} \cup \Aut(\Q.\Z;<)$.

We finally show that the functions generated by $\{\alpha\} \cup \Aut(\Q.\Z;<)$ are dense
in $\Aut(\Q.\Z;\Dist1)$ unless $\alpha \in \Aut(\Q.\Z;F)$. 
When $\alpha$ is not in 
$\Aut(\Q.\Z;F)$, then 
$\alpha$ acts as $i$ on some copies of $\Z$ in $\Q.\Z$,
but not on all. We claim that can assume without loss
of generality that $\alpha$ acts as the identity on infinitely
many copies. Otherwise, $\alpha$ acts as $i$ on infinitely many copies. Then, $\alpha \circ \beta \circ \alpha$
acts as the identity on infinitely many copies,
but as $i$ on some copies. Since $\alpha \circ \beta \circ \alpha$ is generated by $\{\alpha\} \cup \Aut(\Q.\Z;<)$, we can replace $\alpha$ by
$\alpha \circ \beta \circ \alpha$, which satisfies the assumption from the claim. 
Let $\gamma \in \Aut(\Q.\Z;\Dist1)$ be arbitrary. 
Let $q_1,q_2,\dots$ be an enumeration of $\Q$.
We prove by induction on $k \in \N$ that there is a $\alpha'$ generated by $\{\alpha\} \cup \Aut(\Q.\Z;<)$ such that
$\alpha'$ acts as $\gamma$ on 
$q_1.\Z \cup \cdots \cup q_k.\Z$. For $k=1$, 
observe that the restriction of $\gamma$ to $q_1.\Z$
is either in $\Aut(\Q;<)$, or of the form
$\delta \circ i$ for $\delta \in \Aut(\Q;<)$.
In the former case, there is nothing to prove.
In the later case, since there is a copy where $\alpha$ acts as $i$,
we can find the desired $\alpha'$ by composing $\alpha$ with an automorphism of $\Aut(\Q.\Z;<)$. 
In the inductive step, suppose that $\alpha''$ 
is generated by $\{\alpha\} \cup \Aut(\Q.\Z;<)$
and acts as $\gamma$ on $q_1.\Z \cup \cdots \cup q_k.\Z$
First suppose that the restrictions of $\alpha''$ and $\gamma$ 
to $q_{k+1}.\Z$ preserve $\suc$, and that the restrictions
of $\alpha''$ and $\gamma$ to $q_{k+1}.\Z$ both do not preserve $\suc$. In this case we can pick
 $\delta \in \Aut(\Q.\Z;<)$ such that
 $\delta \circ \alpha''$ acts as $\gamma$ on 
 $q_1.\Z \cup \cdots \cup q_{k+1}.\Z$. 
Otherwise, pick $\delta \in \Aut(\Q.\Z;<)$ that maps
$\alpha''(q_1.\Z \cup \cdots \cup q_k.\Z)$ to 
the copies where $\alpha$ behaves as the identity,
and that maps $q_{k+1}.\Z$ to the copy where
$\alpha$ behaves as $i$. Then $\delta^{-1} \circ \alpha \circ \delta \circ \alpha''$ acts as $\Gamma$ on $q_1.\Z \cup \cdots \cup q_k.\Z$,
and preserves $\suc$ on $q_{k+1}.\Z$. By another application of an automorphism of $(\Q.\Z;<)$ we find the desired function generated
by $\{\alpha\} \cup \Aut(\Q.\Z;<)$ that acts as $\gamma$ on $q_1.\Z \cup \cdots \cup q_{k+1}.\Z$. 

This shows in particular that there are no automorphism groups properly between the permutation groups
$\Aut(\Q.\Z;<)$ and $\Aut(\Q.\Z;F)$, and properly between the permutation groups $\Aut(\Q.\Z;F)$ and $\Aut(\Q.\Z; \Dist1)$, finishing the proof. 
\end{proof}

\begin{theorem}\label{thm:rank-one}
Let $\Gamma$ be a finite-signature reduct of $(\Z;<)$ such that $\Q.\Gamma$ has rank one.
Then $\Csp(\Gamma)$ equals $\Csp(\Delta)$
where $\Delta$ is one of the following:
\begin{enumerate}
\item a finite structure;
\item a reduct of $(\Z; <)$ where $\Dist{k}$ is pp-definable for all $k \geq 1$;
\item a reduct of $(\Z;<)$ where $\suc$ is pp-definable.
\end{enumerate}
\end{theorem}
\begin{proof}
If $\Gamma$ has a finite-range endomorphism $f$, then the image of the endomorphism induces a finite structure with the same CSP as $\Gamma$, thus we are in case one and done.
So assume that this is not the case. Then by Theorem~\ref{thm:rank-one-factor}, $\Gamma$ has an endomorphism $g$ that
maps $\Gamma$ to $\Gamma/k$, which is isomorphic to a reduct $\Delta$ of $(\Z;<)$ with the same $\Csp$ as $\Gamma$. 
By Theorem~\ref{thm:rank-one-factor}, all endomorphisms of $\Q.\Delta$ are isometries, and $\Delta$ has no finite-range endomorphism since
by composition of $g$ on $f$ we would obtain a finite-range endomorphism for $\Gamma$.
According to Theorem~\ref{thm:autos}, $\Aut(\Q.\Delta)$ equals 
$\Aut(\Q.\Z;F)$, $\Aut(\Q.\Z;\Dist1)$, or $\Aut(\Q.\Z;<)$. 

Note that since $\Q.\Delta$ has tight-rank one, the relation $\Dist{k}$ is
preserved by the endomorphisms of $\Q.\Delta$.
If $\Aut(\Q.\Delta)$ is either $\Aut(\Q.\Z;F)$ or $\Aut(\Q.\Z;\Dist1)$,
the relation $\Dist{k}$ consists of only one orbit of pairs.
It follows from Theorem~\ref{thm:galois} that $\Dist{k}$ is primitive positive definable in $\Q.\Delta$ for all $k\geq 1$,
and we are in case two of the statement.

In the third case, $\suc$ is preserved by all the endomorphisms of $\Q.\Delta$.
Indeed, let $e$ be an endomorphism of $\Q.\Delta$, and suppose that $e$ does not preserve $\suc$.
Since $\Q.\Delta$ has tight-rank one, we have that $e$ maps copies of $\Z$ to copies of $\Z$.
Composing $e$ with automorphisms of $(\Z;<)$, we may assume that $e(1)-e(0) = -1$,
so that by restricting $e$ to the first copy of $\Z$ in $\Q.\Z$,
we obtain an injective endomorphism $e'$ of $\Delta$ that also violates $\suc$.
By Lemma~\ref{lem:hom-extension}, $e'$ extends to an endomorphism $e''$ of $\Q.\Delta$
 such that $x-y=\infty$ implies $e''(x)-e''(y)=\infty$.
Hence, $e''$ is injective and does not preserve $\suc$.
Note that $e''\circ e''$ acts as a translation on each copy of $\Z$.
Therefore, there exists an automorphism $\alpha$ of $(\Q.\Z;\suc)$
so that $\alpha\circ e''\circ e''$ is the identity on $\Q.\Z$, which means that $e''$
is an embedding.
Finally, the range of $e''$ contains infinitely many
copies of $\Z$ in $\Q.\Z$ by injectivity.
Thus, composing $e''$ with an isomorphism $\iota$ between $e''(\Q.\Delta)$ and $\Q.\Delta$ 
that simply enumerates the copies of $\Z$ that are included in the range of $e''$ by $\Q$, we find that the automorphism $\iota\circ e''$
of $\Q.\Delta$ does not preserve $\suc$, a contradiction to the fact that $\Aut(\Q.\Delta)=\Aut(\Q.\Z;\suc)$.
Since $\suc$ consists of only one orbit of pairs in $\Aut(\Q.\Z;\suc)$, it has a primitive positive definition in $\Delta$ by Theorem~\ref{thm:galois}. 
We are therefore in item three of the statement, which concludes the proof. 
\end{proof}

\subsubsection{Arbitrary Rank}
\label{sect:rank-greater-one}
In this section we study reducts of $(\Z;<)$ with arbitrary finite ranks.

\begin{proposition}\label{prop:rank-k}
Let $\Gamma$ be a reduct of $(\Z;<)$
such that $\Q.\Gamma$ has rank $r \in \N$.
Then $\Gamma/r$ has the same $\Csp$ as $\Gamma$, and is isomorphic to a reduct 
$\Delta$ of $(\Z;<)$ such that 
$\Q.\Delta$ has rank one.
\end{proposition}

Before we can prove the proposition we first have to show some auxiliary results. The following lemma is quite similar,
but formally unrelated, 
to the implication from item~$(\ref{item:saturated-tight-rank})$ to item~$(\ref{item:split})$ in Theorem~\ref{thm:petrus}. 

\begin{lemma}\label{lem:split}
Let $\Gamma$ be a finite-signature reduct of $(\Z;<)$ and $k \in \N$ such that $\Q.\Gamma$ 
is not $k$-bounded. Then for all $x,y \in \Z$ 
such that $x-y = k$ there exists a homomorphism $h$ from $\Gamma$ to $\Q.\Gamma$ 
such that $|h(x) - h(y)| = \infty$. 
\end{lemma}
\begin{proof}
    Since $\Gamma$ is not $k$-bounded, for any $r\geq 0$ there exist $x_0,y_0\in\Z$ and an endomorphism $e\colon\Gamma\to\Gamma$
    such that $e(x_0)-e(y_0)>r$. Composing $e$ with a translation we can take $\{x_0,y_0\}=\{x,y\}$.
    For every finite set $S\subset\Z$, we then have a homomorphism $e\colon\Gamma[S]\to\Gamma$
    such that $e(x)-e(y)>r$.
    Proposition~\ref{prop:koenig} implies   that there exists a homomorphism $h\colon\Gamma\to\Q.\Gamma$
    such that $|h(x)-h(y)|=\infty$.
\end{proof}

\begin{proposition}
\label{prop:rank-image}
Let $\Gamma$ be a finite-signature reduct of $(\Z;\suc)$
such that $\Q.\Gamma$ has rank $r$,
and let $e$ be an endomorphism of $\Q.\Gamma$.
Then $e(z_1) = e(z_2) \mod r$ for all $z_1,z_2 \in \Q.\Z$ such that $z_1=z_2 \mod r$. 
\end{proposition}
\begin{proof}
Suppose that $e \in \End(\Q.\Gamma)$, 
$z_1,z_2 \in \Q.\Z$ contradict the statement of the proposition.
Choose $z_1,z_2$ such that $z_1 > z_2$ and
$z_1-z_2$ is minimal. 

{\bf Claim 1.} $z_1-z_2 = r$.

Suppose otherwise; then there are
$p_1,\dots,p_k$ for $k > 2$ such that
$p_1 = z_1$, $p_k = z_2$, and $p_i - p_{i+1} = r$ for all $i \in \{1,\dots,k-1\}$
because $r$ divides $z_1-z_2$.
By the choice of $z_1,z_2$ we have that
$e(p_i) = e(p_j) \bmod r$. But then 
$e(p_1) = e(p_k) \bmod r$, a contradiction
to the assumption that $e(z_1) \neq e(z_2) \bmod r$. 

Let $w,v \in \mathbb N$ be such that $|e(z_1)-e(z_2)| = wr + v$ and $v < r$. 
Note that $v > 0$ because $e(z_1) \neq e(z_2) \bmod r$.
Assume that $e(z_1) > e(z_2)$; the proof
when $e(z_2) > e(z_1)$ is analogous. 
Let $e' \in \End(\Q.\Gamma)$ be arbitrary, 
and $u_1,u_2 \in \Z$ be arbitrary such that $u_1 - u_2 = v$. 

{\bf Claim 2.} $|e'(u_1) - e'(u_2)| \leq (w+1)c+1$. 

To prove the claim, suppose the contrary.
Let $\alpha  \in \Aut(\Z;<)$ be such that
$\alpha (e(z_1)) = u_1$. Note
that $\alpha(e(z_2) + wr) = u_2$. 
Set $e'' := e' \circ \alpha \circ e$.
Then 
\begin{align*}
& |e''(z_1) - e''(z_2))| \\
\geq \; & |e''(z_1) - e'(u_2)| - |e'(u_2) - e''(z_2))| \\
= \; & |e'(u_1) - e'(u_2)| - |e'(\alpha(e(z_2)+wr)) - e'(\alpha (e(z_2)))| \\
\geq \; & (w+1)c+1 - wc \\
= \; & c+1
\end{align*}
where the first inequality is the triangle inequality,
and the second inequality is by assumption and
$(r,c)$-boundedness.
But $|e''(z_1))-e''(z_2))| > c$ contradicts the assumption that $\Q.\Gamma$ is $(r,c)$-bounded,
and this finishes the proof of Claim 2. 

Since $e'$ was chosen arbitrarily, 
we obtain that $\Q.\Gamma$ is $(v,c)$-bounded,
and hence has rank $v < r$, a contradiction. 
\end{proof}

\begin{lemma}\label{lem:dagger}
Let $\Gamma$ be a finite-signature reduct of $(\Z;<)$ such that $\Q.\Gamma$ has rank $r \in \N$.
Then there exists an endomorphism $e$ of $\Q.\Gamma$ with the following property $(\dagger)$: for all $x \in \Q.\Z$
\begin{align*}
\text{ either } & e(x+k)-e(x) = \infty \\
\text{ or } & e(x+k)-e(x)=0 \mod r
\end{align*}
\end{lemma}
\begin{proof}
We construct $e$ by
an application of K\"onig's tree lemma as follows.
Let $a_1,a_2,\dots$ be an enumeration of the elements of $\Q.\Z$.
The vertices on level $n$ of the tree are $\sim_q$-equivalence classes of homomorphisms $h$ from $\Q.\Gamma[\{a_1,\dots,a_n\}]$ to $\Q.\Gamma$ that can be extended to endomorphisms $h^*$ of $\Q.\Gamma$.
Adjacency between vertices is defined by restriction of representatives. We additionally
require that the homomorphisms $h$ satisfy
property $(\dagger)$ from the statement of the
lemma. 

The interesting part of the proof is to show that the tree has vertices
on all levels. Let $g$ be a homomorphism from $\Q.\Gamma[\{a_1,\dots,a_n\}]$ to $\Q.\Gamma$ that can be extended to an endomorphisms $g^*$ of $\Q.\Gamma$
such that the number $m$ of pairs $i,j \in \{1,\dots,n\}$
with $g(a_i) - g(a_j) = \infty$ or 
$g(a_i) = g(a_j) \mod r$ is maximal. 
If $m = {n \choose 2}$ then we are done; so suppose that there are $p,q \in \{1,\dots,n\}$
such that $g(a_p)-g(a_q) \in \Z$ is not divisible by $r$. Let $k \in \{1,\dots,r-1\}$ and $l \in \Z$ be such that $g(a_p)-g(a_q) = lr+k$, $0<k<r$. 
Since $\Q.\Gamma$ is not $k$-bounded, 
by Lemma~\ref{lem:split}
there exists a homomorphism $f'$ from the copy
of $\Gamma$ that contains $g(a_p)$, $g(a_p)+lr$, and $g(a_q)$ to $\Q.\Gamma$
such that $f'(g(a_p)+lr) - f'(g(a_q)) = \infty$.
We assume without loss of generality that the
range of $f'$ is bounded, in the sense that the set of copies of $\Z$ that intersect the image of $f'$ is bounded. 
By Lemma~\ref{lem:hom-extension} there exists an endomorphism $f$ of $\Q.\Gamma$ that extends $f'$ and has the property
that $f(x) - f(y) = \infty$ whenever $x-y = \infty$. 

By Proposition~\ref{prop:rank-image} we have that $f(g(a_p)) = f(g(a_p)+lr) \mod r$, 
and hence $f(g(a_p)) - f(g(a_q)) = \infty$.
We claim that the number $m'$ of pairs $i,j \in \{1,\dots,n\}$
such that $f(g(a_i)) - f(g(a_j)) = \infty$ or 
$f(g(a_i)) = f(g(a_j)) \mod r$ is larger than $m$.
If $g(a_i) - g(a_j) = \infty$ then $f(g(a_i))-f(g(a_j)) = \infty$; if $g(a_i) = g(a_j) \mod r$ then $f(g(a_i)) = f(g(a_j)) \mod r$. Therefore, $m' \geq m$. Moreover, we have $f(g(a_p)) - f(g(a_q)) = \infty$, and hence $m' > m$.
But $f \circ g^*$ is an endomorphism of $\Q.\Gamma$, contradicting the maximality of $m$.
\end{proof}

We will now prove Proposition~\ref{prop:rank-k}.
Let $\Gamma$ be a finite-signature reduct of $(\Z;<)$
such that $\Q.\Gamma$ has rank $r \in \N$.
We will show that $\Gamma/r$ has the same CSP as $\Gamma$, and is isomorphic to a reduct 
$\Delta$ of $(\Z;<)$ such that 
$\Q.\Delta$ has rank one.

\begin{proof}
By Lemma~\ref{lem:division}, 
there is a reduct $\Delta$ of $(\Z;<)$
such that 
$x \mapsto r \cdot x$ is an isomorphism
between $\Delta$ and $\Gamma/r$. 
Let $e$ be the endomorphism of $\Q.\Gamma$ constructed in 
Lemma~\ref{lem:dagger}. Replacing 
$e$ by $\alpha \circ e$ for an appropriate
automorphism of $(\Q.\Z;<)$, we
can assume that the range of $e$ lies within
$S:=\{r \cdot z : z \in \Q.\Z\}$. 
Since $x \mapsto r \cdot x$ is an isomorphism between $\Q.\Delta$
and the structure
induced by $S$ in $\Q.\Gamma$,
we obtain that $\Gamma$, $\Q.\Gamma$,
$\Q.\Delta$, and $\Delta$ all have
the same $\Csp$. 

It remains to be shown that 
$\Q.\Delta$ has rank $1$. 
Suppose for contradiction that $\Q.\Delta$
is not 1-bounded. 
Then by Lemma~\ref{lem:split} there exists
an $f \in \End(\Q.\Delta)$ and an $z \in \Q.\Z$ such that $f(z+1)-f(z) = \infty$.
By composing $f$ with an endomorphism of $(\Q.\Z;<)$, we can assume that the range of $f$ is bounded. 
The map $f' \colon S \to \Q.\Z$ given by $x \mapsto f(x/r)$ is a homomorphism
from $\Q.\Gamma[S]$ to $\Q.\Gamma$, which can be extended to an endomorphism
$g$ of $\Q.\Gamma$ by Lemma~\ref{lem:hom-extension}. Then $$g(r \cdot z + r) - g(r \cdot z) = f(z+1) - f(z) = \infty \; ,$$
and hence $\Q.\Gamma$ is not $r$-bounded, in contradiction to our rank $r$ assumption on $\Q.\Gamma$. 
\end{proof}

\subsection{Defining $\suc$ and $<$}

Combining the results of the preceding subsections, we finally get the following:

\begin{thm:class}
Let $\Gamma$ be a reduct of $(\Z;<)$ with finite signature. 
Then $\Csp(\Gamma)$ equals $\Csp(\Delta)$
where $\Delta$ is one of the following:
\begin{enumerate}
\item a finite structure;
\item a reduct of $(\Q;<)$;
\item a reduct of $(\Z;<)$ where $\Dist{k}$ is pp-definable for all $k \geq 1$;
\item a reduct of $(\Z;<)$ where $\suc$ is pp-definable. 
\end{enumerate}
\end{thm:class}
\begin{proof}
Let $\Gamma$ be a reduct of $(\Z;<)$ with finite signature. If $\Gamma$ has an endomorphism with finite range, then
$\Gamma$ is homomorphically equivalent to a finite structure; hence item~1 of Theorem~\ref{thm:class} applies and we are done.
So suppose that this is not the case. If there exists a reduct of $(\Q;<)$ with the same CSP, then item~2 of Theorem~\ref{thm:class} applies and we are done.
Otherwise, the equivalence of $(\ref{item:saturated-tight-rank})$ and $(\ref{item:equality})$ in Theorem~\ref{thm:petrus} implies that $\Q.\Gamma$
has bounded tight rank $t$ and bounded rank $r$. If $r = 1$, then the statement follows from Theorem~\ref{thm:rank-one}. 
Otherwise, if $r > 1$, then by Proposition~\ref{prop:rank-k}, $\Gamma$ has the same $\Csp$ as $\Gamma/r$ and $\Gamma/r$ is isomorphic to a reduct $\Delta$ of $(\Z;<)$ such that $\Q.\Delta$ has rank 1.
The statement again follows from Theorem~\ref{thm:rank-one}.
\end{proof}

In the rest of this section, we now prove the following dichotomy:
a reduct of $(\Z;<)$ that pp-defines $\suc$ either pp-defines $<$, or is a first-order
reduct of $(\Z;\suc)$.
Call a binary relation $R$, fo-definable in $(\Z;<)$, \emph{one-sided infinite} if there exists $c\leq d$ so that
for all $z< c$, $R(x,x+z)$ does not hold, but, for all $z\geq d$, $R(x,x+z)$ holds.
The following is clear.
\begin{lemma}\label{lem:one-sided-infinite}
For a reduct $\Gamma$ of $(\mathbb{Q}.\Z;<)$ in which $\suc$ is pp-definable, $<$ is pp-definable iff some one-sided infinite binary relation is pp-definable.
\end{lemma}
\begin{proof}
Since $<$ is one-sided infinite we only have to show the reverse implication. Let $R$ be a binary one-sided infinite relation with a pp-definition in $\Gamma$,
and choose $c,d$ as in the definition such that $d-c$ is minimal.
If $c=d$ then $R$ is a relation of the form $x<y+k$ for $k\in\Z$, and using $\suc$ we can pp-define $<$ in $\Gamma$.
So assume that $c\neq d$.
Replace $R$ by the relation $T$ defined by the formula $R(x,y)\land R(x,y+d-c-1)$, which is equivalent to a pp-formula over $\Gamma$.
Then $(0,x)$ is in $T$ for all $x\geq d$. On the other hand, for $x<c+1$, we have that $(0,x)\not\in T$.
Indeed, if $x<c$ then $(0,x)$ is not in $R$, so not in $T$. If $x=c$, then $(0,d-1)$ is not in $R$ by the minimality of $d$, so that $(0,c)$ is not in $T$.
Therefore the integers $c,d$, as defined for $T$, have a smaller difference than those for $R$.
We can then proceed by induction until $c=d$.
\end{proof}

If $R$ is a relation of arity $n$, and $i_1,\dots,i_k\in\{1,\dots,n\}$ are distinct indices,
the projection of $R$ onto $\{i_1,\dots,i_k\}$, denoted by $\pi_{i_1,\dots,i_k}(R)$,
is the relation defined by $\exists_{j\not\in\{i_1,\dots,i_k\}} x_j. R(x_1,\dots,x_n)$ over $(\Z;R)$.
A \emph{binary projection} of $R$ is a projection of $R$ onto a set of size $2$.

\begin{lemma} 
\label{lem:define-order}
Let $\Gamma$ be a finite-signature reduct of $(\Z;<)$ in which $\suc$ is pp-definable. Then, either $\Gamma$ pp-defines $<$ or $\Gamma$ is a reduct of $(\Z;\suc)$.
\end{lemma}
\begin{proof}
Take a relation $R(x_1,x_2,x_3,\ldots,x_k)$ of $\Gamma$.
If $E=\{x-y \mid (x,y)\in\pi_{i,j}(R)\}$ is a finite or cofinite set,
we can define $R$ over $(\Z;<)$ without ever using a literal of the form $x_i<x_j+k$.
Indeed, such a literal can be replaced by a disjunction of literals $\suc^p(x_i,x_j)$ for suitable
integers $p$ if $E$ is finite, or a by a conjunction of literals $\neg\suc^p(x_i,x_j)$
when $E$ is cofinite. Therefore if $\Gamma$ is not a reduct of $(\Z;\suc)$
there exists a relation $R$ of $\Gamma$ and integers $i,j$ such that $\pi_{i,j}(R)$
is a binary relation and such that the set $\{x-y\mid (x,y)\in\pi_{i,j}(R)\}$ is neither finite nor cofinite.
Since $R$ has a finite qe-degree, it follows that $R$ or $\{(y,x) \in\Z^2\mid (x,y)\in R\}$ is one-sided infinite.
From Lemma~\ref{lem:one-sided-infinite} and the fact that $\pi_{i,j}(R)$ is pp-definable in $\Gamma$ follows that $<$ is pp-definable in $\Gamma$.
\end{proof}

\begin{corollary}
Let $\Gamma$ be a finite-signature reduct of $(\Z;<)$. There exists a structure $\Delta$
such that $\Csp(\Delta)=\Csp(\Gamma)$ and such that one among the following holds.
\begin{itemize}
\item $\Delta$ is a finite structure,
\item $\Delta$ is a reduct of $(\Q;<)$,
\item $\Delta$ is a reduct of $(\Z;\suc)$ that pp-defines $\suc$ or $\Dist{k}$ for all $k\geq 1$.
\item $\Delta$ is a reduct of $(\Z;<)$ that pp-defines both $\suc$ and $<$.
\end{itemize}
\end{corollary}
\begin{proof}
Follows from culmination of previous section (Theorem~\ref{thm:class}) and Lemma~\ref{lem:define-order}.
\end{proof}

\section{Tractable classes}
\label{sect:algs}
We treat items 3 and 4 in Theorem~\ref{thm:main}, that is, we prove that if $\Gamma$ is a reduct of $(\Z;<)$ that is preserved by $\max_d$ or $\min_d$,
or if $\Gamma$ is a reduct of $(\Z;\suc)$ such that $\Q.\Gamma$ is preserved by a binary injective operation preserving $\suc$,
then $\Csp(\Gamma)$ is in Ptime.

\subsection{The Horn Case}
Let $\si$ be any isomorphism between $(\Q.\Z,\suc)^2$ and $(\Q.\Z,\suc)$.
Remember that relations that are first-order definable in $(\Z;\suc)$
are also definable
by quantifier-free formulas with (positive or negative) literals of the form
$\suc^p(x,y)$ for $p\in\Z$ (see Example~\ref{ex:F}).
A quantifier-free formula in conjunctive normal form over $\suc$ is called \emph{Horn}
 if each clause of the formula contains at most one positive literal,
that is, at most one literal of the form $\suc^p(x,y)$.
A relation is said to be \emph{Horn-definable in $(\Z;\suc)$} if there exists a Horn formula that defines
the relation in $(\Z;\suc)$.

\begin{proposition}\label{prop:si_preserved_implies_horn}
Let $\Gamma$ be a reduct of $(\Z;\suc)$. 
If $\Q.\Gamma$ is preserved by
$\si$ then every relation of $\Gamma$ has a quantifier-free Horn definition over
$(\Z;\suc)$.
\end{proposition}
\begin{proof}
Since $(\Z;\suc)$ has quantifier elimination, the result follows from Proposition~5.9 of \cite{BodHilsMartin-Journal}. 
\end{proof}

\begin{proposition}\label{prop:horn_implies_tractable}
Let $\Gamma$ be a reduct of $(\Z;\suc)$ such that every relation of $\Gamma$ has a quantifier-free Horn definition over $(\Z;\suc)$. Then $\Csp(\Gamma)$ is in Ptime. 
\end{proposition}
\begin{proof}
It is easy to see that there is an algorithm
that decides whether a set of constraints of the form
$\suc^{p_i}(x_i,y_i)$ is consistent or if it implies another constraint of this form.
Indeed, to see if the set of constraints is consistent,
consider the graph whose vertices are the variables,
and whose arcs are $(x_i,y_i)$, labelled by $p_i$,
if there is a constraint $\suc^{p_i}(x_i,y_i)$.
For each variable $x$, using a graph traversal we can
check if all the directed paths going from $x$
to some other variable $y$ have the same weight (which is given by the sum of the labels over the arcs);
If this is not the case, the constraints are unsatisfiable. 
Otherwise, to decide whether the constraints imply $\suc^p(x,y)$,
check if there is a directed path from $x$ to $y$ where the sum of the labels equals $p$.

We view the instance of $\Csp(\Gamma)$ as a 
set of Horn-clauses over $(\Z;\suc)$, and run positive unit resolution on this set, using the above
algorithm to test whether a literal can be eliminated from a clause. If we derive
the empty clause, reject the input. Otherwise, 
the resolution stabilizes in a polynomial number of steps with a set of Horn clauses; in this case,
 accept the input. 
We apply $\si$ to show that in this case indeed there exists a solution.
By assumption, for each Horn clause $\bigwedge \suc^{p_i}(x_i,y_i) \Rightarrow \suc^p(x,y)$ there exists an  assignment that falsifies some literal $\suc^{p_i}(x_i,y_i)$ and additionally satisfies all the positive unit clauses: otherwise the literal would have been removed by the resolution procedure.
Let $s_1,\ldots,s_r$ be those assignments for 
the $r$ clauses.
Since $\si$ is an isomorphism,
the assignment $s := \si(s_1,\ldots,\si(s_{r-1},s_r)\ldots)$ simultaneously breaks 
all the equalities in the premises of all the clauses.
Moreover, since $\si$ preserves $\suc$, 
the resulting assignment $s$ also preserves the
positive unit clauses, and hence is a valid assignment for the input. 
\end{proof}

Putting together Propositions~\ref{prop:si_preserved_implies_horn} and \ref{prop:horn_implies_tractable},
we obtain the following corollary:

\begin{corollary}\label{cor:si-tract}
	Let $\Gamma$ be a reduct of $(\Z;\suc)$.
	If $\Q.\Gamma$ is preserved by $\si$, then $\Csp(\Gamma)$ is in Ptime.
\end{corollary}
 
\subsection{Modular Minimum and Modular Maximum}

The \emph{modular minimum} and \emph{modular maximum} were defined in~\cite{BodDalMarMotPin}.

\begin{definition}
		Let $d$ be a positive integer. 
		The \emph{$d$-modular max} is the binary operation $\max_d\colon\Z^2\to\Z$
		that is defined by $\max_d(x,y) := \max(x,y)$ if $x=y\mod d$ and $\max_d(x,y):=x$ otherwise.
		The \emph{$d$-modular min} is similarly defined as the operation $\min_d\colon\Z^2\to\Z$
		which satisfies $\min_d(x,y):=\min(x,y)$ if $x=y\mod d$ and $\min_d(x,y):=x$ otherwise.
\end{definition}

\begin{theorem}\label{thm:tractability-modular-semilattice}
	Let $\Gamma$ be a finite-signature reduct of $(\Z;<)$
	that admits a modular max or modular min polymorphism.
	Then $\Csp(\Gamma)$ is in Ptime.
\end{theorem}
\begin{proof}
	Suppose that $\Gamma$ is preserved by $\max$, the regular maximum operation.
	Then $\Csp(\Gamma)$ is solvable in polynomial-time as follows.
	Let $q$ be the qe-degree of $\Gamma$.
	Let $\phi$ be an instance of $\Csp(\Gamma)$ with $n$ variables.
	We already noted in the proof of Proposition~\ref{prop:inNP}
	that $\phi$ is satisfiable in $\Gamma$
	iff it is satisfiable in $\Gamma[\{0,\dots,(q+1)n\}]$,
	and the latter structure can be constructed in polynomial-time,
	and is preserved by the maximum function on $\{0,\dots,(q+1)n\}$.
	We can then decide whether $\Gamma[\{0,\dots,(q+1)n\}]\models\phi$
	using the arc-consistency algorithm, noting that the arc-consistency procedure
    can be implemented in such a way that the running time is linear in both the size of the formula and of the structure.
	
	Suppose now that $\Gamma$ is preserved by $\max_d$ for $d>2$.
	It follows that $<$ is not pp-definable in $\Gamma$, as $\max_d$ does not preserve $<$.
	We can suppose that $\Gamma$ pp-defines $\suc$, because
	this only increases the complexity of $\Csp(\Gamma)$ and $\suc$ is preserved by $\max_d$.
	By Lemma~\ref{lem:define-order}, $\Gamma$ is a first-order reduct of $(\Z;\suc)$.
	In~\cite{BodDalMarMotPin}, the authors prove
	that the CSP of a first-order reduct of $(\Z;\suc)$ with finite distance degree
and which is preserved by a modular maximum or minimum is decidable in polynomial-time.
An inspection of the proof shows that the finite distance degree hypothesis is not necessary.
Indeed, the critical idea of the algorithm is that if $\Gamma$ is preserved by the $d$-modular maximum, then $\Csp(\Gamma)$ reduces in polynomial time to $\Csp(\Delta)$,
where $\Delta$ is a reduct of $(\Z;\suc)$ which is preserved by the usual maximum or minimum. This reduction does not rely on the distance degree of $\Gamma$ being
finite to work.
\end{proof}

\section{The Classification}
\label{sect:classification}
In this section we prove Theorem~\ref{thm:main}. By Theorem~\ref{thm:class},
we are essentially left with the task to 
classify the CSP for finite-signature reducts $\Gamma$ of $(\Z;<)$ where the binary relation $\suc$ is among the relations of $\Gamma$.

\begin{lemma}\label{lem:not-max-NPhard}
	Let $\Gamma$ be a reduct of $(\Z;<)$ that pp-defines $\suc$ and $<$.
	If\/ $\Gamma$ is preserved by neither $\max$ nor $\min$, then $\Csp(\Gamma)$
	is NP-hard.
\end{lemma}
\begin{proof}
	Let $R,T$ be a pair of relations of $\Gamma$ which are not preserved
	by $\max$ and $\min$, respectively.
	Therefore there are tuples $\tuple a,\tuple b$ in $R$ (resp. $T$)
	such that $\max(\tuple a,\tuple b)\not\in R$ (resp. $\min(\tuple a,\tuple b)\not\in T$).
	Let $M$ be $\max_{i,j}(|a_i-a_j|,|b_i-b_j|)$.
    Since $\suc$ and $<$ have pp-definitions in $\Gamma$,
    the binary relation defined by $x\leq y+M$ has a pp-definition in $\Gamma$ as well.
	We define $R^*$ by the following  primitive positive formula in $\Gamma$
	$$R(x_1,\dots,x_n)\land \bigwedge_{i,j} x_i \leq x_j + M$$
	and note that $\tuple a$ and $\tuple b$ are in $R^*$, and that $\max(\tuple a,\tuple b)$
	is not in $R^*$. Note that $R^*$ is first-order definable over $\suc$
	and has finite distance degree.
	Similarly we define a relation $T^*$ which is pp-definable in $\Gamma$,
	not preserved by $\min$, and which is first-order definable over $\suc$.
	It follows from Proposition 47 in~\cite{BodDalMarMotPin}
	that $\Csp(\Z;\suc,R^*,T^*)$ is NP-hard, so that $\Csp(\Gamma)$
	is NP-hard.
\end{proof}

By Lemma~\ref{lem:define-order}, if $\Gamma$ is a first-order reduct of $(\Z;<)$
which defines $\suc$ but not $<$, then $\Gamma$ is a reduct of $(\Z;\suc)$.
We call such a structure $\Gamma$ a \emph{first-order (fo-) expansion} of $(\Z;\suc)$. 
In the following, we will use the relations of the form $\suc^p$ (see Example~\ref{ex:F}) as if they were atomic symbols of the language.
Since they are all pp-definable in an fo-expansion of $(\Z;\suc)$, this will not cause any loss of generality.

\begin{theorem}\label{thm:succ}
Let $\Gamma$ be a first-order expansion of $(\Z;\suc)$. 
Then at least one of the following is true:
\begin{enumerate}
	\item $\Gamma$ is preserved by a modular max or a modular min,
	\item $\Q.\Gamma$ is preserved by $\si$,
	\item $\Csp(\Gamma)$ is NP-hard.
\end{enumerate}
\end{theorem}

Binary relations $R$ with a first-order definition in $(\Z;\suc)$
come in two flavours. Indeed, the set $\{x-y \mid (x,y)\in R\}$ is either \emph{finite} or \emph{cofinite} and by abuse of language
we will say that $R$ is finite or cofinite.
A binary relation $R\subseteq\Z^2$ that is first-order definable in $(\Z;\suc)$
is called \emph{trivial} if it is pp-definable over $(\Z;\suc)$, and \emph{non-trivial} otherwise. 
Moreover, we also call a binary relation that is first-order definable over $(\Q.\Z;<)$
trivial if it is pp-definable over $(\Q.\Z;\suc)$, and non-trivial otherwise.
We first give a syntactic characterization of those reducts of $(\Z;\suc)$
in which no non-trivial binary cofinite relation is pp-definable.
A formula over $\suc$ is said to be \emph{positive}
if it only includes positive literals of the form $\suc^p(x,y)$.
A formula over the signature of $(\Z;\suc)$ in DNF is called \emph{reduced} when 
every formula obtained by removing literals or conjunctive clauses is not logically equivalent over $(\Z;\suc)$.
It is clear that every first-order formula on $(\Z;\suc)$ is equivalent to a reduced formula in DNF.

\begin{lemma}\label{lem:positive-iff-endo-collapses}
	For a first-order expansion $\Gamma$ of $(\Q.\Z;\suc)$,
	the following are equivalent:
	\begin{enumerate}
		\item\label{itm:reduced-DNF-positive} Every reduced DNF that defines a relation of $\Gamma$ is positive;
		\item\label{itm:collapsing-endo} $\Gamma$ has an endomorphism that violates the binary relation given by $\lvert x-y\rvert = \infty$;
		\item\label{itm:no-binary-cofinite-non-trivial} $\Gamma$ does not pp-define a non-trivial binary relation with infinite distance degree.
	\end{enumerate}
\end{lemma}
\begin{proof}
    $(\ref{itm:reduced-DNF-positive})$ implies $(\ref{itm:collapsing-endo})$. Any function preserving $\suc$
    is an endomorphism of $\Gamma$, in particular $\Gamma$ admits an endomorphism violating $|x-y|=\infty$.

	(\ref{itm:collapsing-endo}) implies (\ref{itm:reduced-DNF-positive}).
	Let $e$ be an endomorphism of $\Gamma$ that violates $x-y=\infty$, and let
	$a,b$ be such that $a-b=\infty$ and $e(a) - e(b) <\infty$.
	Using automorphisms of $(\Q.\Z;\suc)$, we may assume that $e(a)=e(b)=b$ without loss of generality.
	For contradiction, suppose that $\Gamma$
	has a relation with a reduced DNF definition $\phi(x_1,\ldots,x_n)$ which is not positive.
		
	We now show that we can choose $s\colon\{x_1,\dots,x_n\}\to\Z$ such that $s$ is a satisfying assignment for $\phi$
	but $e\circ s$ is not.
	For this, let us write one of the non-positive disjuncts $\psi$ of $\phi$ as $\neg\suc^p(z_2,z_1)\land\phi'$
	where $\phi'$ is a conjunction of literals and $z_1,z_2\in\{x_1,\dots,x_n\}$.
	Moreover, let $\psi_2,\ldots,\psi_m$ be the other disjuncts of $\phi$.
	Suppose that all assignments that satisfy $\phi'\land\suc^p(z_2,z_1)$
	also satisfy $\bigvee_{2\leq i\leq m} \psi_i$.
	Then we could rewrite $\phi$ as simply $\phi'\lor\bigvee \psi_i$,
	which is impossible since $\phi$ is reduced.
	Hence, there exists $t\colon\{x_1,\dots,x_n\}\to\Z$
	such that $t$ is a satisfying assignment for
	$\phi'\land \suc^p(z_2,z_1)$ but is not a satisfying assignment for any $\psi_i$.
	Using an automorphism of $(\Q.\Z;\suc)$, we can assume that $t(z_2) = b-p$.
	Moreover, we can assume that $t$ occupies only one copy of $\Z$:
	let $S$ be the set $\{t(x_1),\dots,t(x_n)\}$ and let $g\colon S\rightarrow\Q.\Z$
	be any function that maps $S$ to the first copy of $\Z$ in such a way that if $t(x_i)$ and $t(x_j)$
	are in different copies, then $g(t(x_i))$ and $g(t(x_j))$ are at distance at least $q+1$,
	where $q$ is the qe-degree of $\phi$.
	We have that $g$ is $\sim_q$-equivalent to any embedding of $S$ into the first copy of $\Z$ in $\Q.\Z$.
	Therefore by the Substitution Lemma (Lemma~\ref{lem:substitution}), the function $g\circ t$ is a satisfying assignment to the variables of $\phi$
	that only occupies one copy of $\Z$.
	
	We now derive from $t$ an assignment $s$ that satisfies $\neg\suc^p(z_2,z_1)$,
	that agrees with $t$ on all the other literals of $\psi$
	and such that $e(s) = t$:
	if we consider $\phi'$ as a graph on $\{z_1,\ldots,z_k\}$ where edges represent positive literals,
	$z_1$ and $z_2$ are in different connected components. Indeed, if there were a path from $z_1$ to $z_2$
	in this graph we would have that $\phi'$ implies a statement of the form
	$\suc^q(z_2,z_1)$. But then
	the conjunction $\neg\suc^p(z_2,z_1)\land\suc^q(z_2,z_1)$ is either contradictory or is equivalent
	to $\suc^q(z_2,z_1)$, which is a contradiction since $\phi$ is reduced.
	Let $V$ be the variables in the connected component of $z_1$,
	and define $s$ on $V$ as $s(v) = a - t(z_1) + t(v)$ (in particular $s(z_1) = a$)
	and $s(v) = t(v)$ on the variables that are not in $V$.
	We have that $s$ satisfies $\neg\suc^p(z_2,z_1)$ and that $s$
	agrees with $t$ on the other literals in $\psi$: the truth of positive literals
	is preserved since we performed a translation on variables that are connected by positive literals,
	and negative literals between the variables in $V$ and the other variables are trivially true,
	since the distance between a variable in $V$ and a variable not in $V$ is infinite.
	Finally, a negative literal between variables not in $V$ is preserved since those variables are not touched.
	Hence, $s$ is a satisfying assignment of $\phi$.
	We have $e\circ s = t$.
	If $v$ is a variable in $V$, then $e(s(v)) = e(a) - t(z_1) + t(v) = t(v)$, and if $v\not\in V$ we \emph{defined} $s(v)$ to be $t(v)$,
	so that $e(s(v)) = e(t(v)) = t(v)$.
	This contradicts the fact that $e$ is an endomorphism of $\Gamma$.
	
    (\ref{itm:reduced-DNF-positive}) implies (\ref{itm:no-binary-cofinite-non-trivial}). Let $R$ be a binary relation with a pp definition in $\Gamma$
	of the form $\exists \overline z \bigwedge_i R_i(x,y,\overline z)$.
	Let us replace $R_i(x,y,\overline z)$ by a quantifier-free formula $\phi_i$ in reduced DNF that defines $R_i$ in $(\Q.\Z;\suc)$.
	By assumption, all the literals in $\phi_i$ are positive.
	This formula is equivalent to $\psi(x,y) := \bigvee_i \exists\overline z \psi_i(x,y,\overline z)$ where $\psi_i$ is a conjunction of atoms.
	Note that one can eliminate quantifiers in $(\Q.\Z;\suc)$ in such a way that pp formulas are translated into positive qf formulas.
	If one of the disjuncts in $\psi(x,y)$ is vacuously true after quantifier elimination, then $\psi$ defines a trivial binary relation.
	Otherwise, all the literals are $\suc^{p_i}(y,x)$ for suitable $p_i\in\Z$, and the relation defined by $\psi$ is finite.
	In either case, $\psi$ does not define a non-trivial binary relation of infinite distance degree.

	(\ref{itm:no-binary-cofinite-non-trivial}) implies (\ref{itm:reduced-DNF-positive}).
    Suppose that $\Gamma$ has a relation $R$ whose reduced DNF definition $\phi$ over $(\Z;\suc)$
    is not positive. Assume that $\phi$ contains a minimal number of negative literals.
    There exists a clause of $\phi$ which contains a negative literal, say $\neg\suc^p(x_j,x_i)$.
    Let $C$ be the conjunction of the other literals in this clause.
    Write $\phi'$ for the rest of the formula.
    Since $\phi$ is reduced, there exists a tuple $\tuple a$ in $\Z^n$ that satisfies $\neg\phi'\land C\land \suc^p(x_j,x_i)$
    (otherwise $\phi$ is equivalent to $\phi'\lor C$).
    As above we consider $C$ as a graph on the variables, and we again have that there is no path from $x_j$ to $x_i$.
    Let $V_1$ be the connected component of $x_i$ and $V_2$ be the set of all the other variables.
    Let $P$ be the conjunction of literals of the form $\suc^p(x_i,x_k)$ where $a_k=a_i+p$ with $x_k\in V_1$,
    and of the form $\suc^p(x_j,x_k)$ where $a_k=a_j+p$ with $x_k\in V_2$.
    Let $\psi(x_i,x_j)$ be the formula $\exists_{k\not\in\{i,j\}} x_k (R(x_1,\dots,x_n)\land P(x_1,\dots,x_n))$.
    We have that $\psi$ is not satisfied by $(a_i,a_j)$: the existentially quantified variable $x_k$
    would need to be instanciated by $a_k$ for every $k\not\in\{i,j\}$ but we know that $\tuple a$ is not in $R$.
    Moreover $(a_i,a_i+M)$ satisfies $\psi$ for $M$ large enough: instanciate the variable $x_k$ in $V_1$ by $b_k:=a_k$ and the variable $x_k$ in $V_2$ by $b_k:=a_k+M$.
    This new tuple $\tuple b$ satisfies all the literals in $P$, but also all the literals in $C$: a literal
    whose variables are both in $V_1$ or both in $V_2$ is satisfied by $\tuple b$, and a literal whose variables
    are in different $V$-sets is necessarily negative by construction.
    Such a literal is satisfied provided that $M$ is large enough.
    Therefore the relation defined by $\psi$ has infinite distance degree and is not $\Z^2$.
    This concludes the proof of the lemma.
\end{proof}

Using Lemma~\ref{lem:positive-iff-endo-collapses},
we treat positive and non-positive
expansions of $(\Z;\suc)$ separately in subsections 7.1 and 7.2.

\subsection{The non-positive case}

Our aim in this section is to show that a non-positive reduct $\Gamma$ of $(\Z;\suc)$ has an NP-hard CSP
if $\Q.\Gamma$ is not preserved by $\si$.
In order to do this, we show that there is a non-trivial binary relation with finite distance degree
that is pp-definable in $\Gamma$.

We work in the following with reduced standardized formulas. We say that a formula $\phi$ in CNF is \emph{reduced}
when removing any literal in a clause yields a formula that is not equivalent to $\phi$.
This is equivalent to say that for any literal $\ell$ in a clause $\psi$ of $\phi$, there exists an assignment that satisfies $\phi$
and that satisfies only $\ell$ in $\psi$. This assignment witnesses the fact that the given
literal cannot be removed from the formula without changing the set of satisfying assignments.
Given two formulas $\phi$ and $\psi$, we say that we obtain $\psi$ by reducing $\phi$
if $\psi$ is obtained from $\phi$ by removing clauses or literals and if $\phi$ and $\psi$ are equivalent.

\begin{lemma}\label{lem:Horn-iff-reduced-Horn}
	Let $\phi$ be a formula in the language of $(\Z;\suc)$, and
	suppose that $\phi$ is equivalent to a Horn formula over $(\Z;\suc)$.
	Then any reduced formula obtained by reducing $\phi$ is Horn.
\end{lemma}
\begin{proof}
	Note that $\phi$ is equivalent to a Horn formula over $(\Z;\suc)$
	if and only if it is equivalent to a Horn formula over $(\Q.\Z;\suc)$,
	since both structures have the same first-order theory.
	We therefore have that $\phi$ is preserved by $\si$.
	Let $\psi$ be a reduced standardized formula obtained by reducing $\phi$.
	Suppose for contradiction that $\psi$ is not Horn, that is, it contains a clause $\psi'$ of the form
	$$(\suc^p(y,x) \lor \suc^q(v,u) \lor \dots).$$
	Since this formula is reduced, there exist assignments $s,t$ that satisfy $\psi$
	and such that $s$ satisfies only $\suc^p(y,x)$ in $\psi'$,
	and $t$ only satisfies $\suc^q(v,u)$
	in $\psi'$.
	The assignment $(s,t)$ that maps a variable $x_i$ of $\psi$ to the pair $(s(x_i),t(x_i))$ in $(\Q.\Z)^2$
	is not a satisfying assignment for $\psi$. Since $\si$ is an isomorphism between $(\Q.\Z;\suc)^2$ and $(\Q.\Z;\suc)$,
	we have that the assignment $\si(s,t)$ does not satisfy $\psi$, which contradicts the fact that $\psi$
	is preserved by $\si$.
\end{proof}

\begin{proposition}\label{prop:non-Horn-binary-finite}
	Let $\Gamma$ be a first-order expansion of $(\Z;\suc)$,
	and suppose that $\Gamma$ pp-defines a relation that is not Horn-definable over $(\Z;\suc)$.
	Then $\Gamma$ also pp-defines a binary relation that is not Horn-definable over $(\Z;\suc)$.
\end{proposition}
\begin{proof}
	Let $R$ be a relation with a pp-definition in $\Gamma$ that is not Horn-definable over $(\Z;\suc)$,
	and whose arity is minimal among the relations with same properties.
	We claim that $R$ is binary.
	Let $\phi(x_1,\dots,x_n)$ be a reduced standardized formula that defines $R$ in $(\Z;\suc)$ whose number of non-Horn clauses is minimal,
	and suppose for contradiction that $n>2$.
	
	We first prove that $\phi$ has a non-Horn clause that consists of positive literals only.
	Pick a non-Horn clause $\psi$ of $\phi$, and suppose it contains the negative literal $\neg\suc^p(x_j,x_i)$
	for some $i,j\in\{1,\dots,n\}$ and $p\in\Z$.
	We claim that the formula $\phi\land \suc^p(x_j,x_i)$ is not equivalent to a Horn formula.
	By Lemma~\ref{lem:Horn-iff-reduced-Horn}, it suffices to prove that when reducing $\phi\land \suc^p(x_j,x_i)$,
	there remains at least two positive literals in $\psi$.
	Pick two positive literals $\ell_1$ and $\ell_2$ in $\psi$. Since $\phi$ is reduced, there is an assignment $s_1$ (resp. $s_2$)
	that satisfies $\phi$, and that satisfies only $\ell_1$ (resp. $\ell_2$) in $\psi$.
	In particular, this assignment does not satisfy $\neg\suc^p(x_j,x_i)$, so that it satisfies $\suc^p(x_j,x_i)$.
	As a consequence, both $s_1$ and $s_2$ are satisfying assignments of
	$\phi\land \suc^p(x_j,x_i)$
	that satisfy exactly one literal in $\psi$, which means that $\ell_1$ and $\ell_2$
	cannot be removed when reducing $\phi\land \suc^p(x_j,x_i)$.
	
	From the previous paragraph we get that $\phi\land \suc^p(x_j,x_i)$ is not equivalent to a Horn formula.
	Using the same argument, we see that $\exists x_j(\phi\land \suc^p(x_j,x_i))$
	is not equivalent to a Horn formula either. Note that this formula defines a relation of arity $n-1$ that is not Horn
	and that is pp-definable in $\Gamma$, a contradiction to the choice of $R$.

	Therefore, there exists a positive clause $\psi$ in $\phi$. We let $\phi'$ denote the rest of the formula.
	Let $E_{i,j}$ be the set
	$$\{s(x_j)-s(x_i) \mid s\colon\{x_1,\dots,x_n\}\to\Z \text{ satisfies } \phi'\land\neg\psi\}.$$
	If $E_{i,j}$ is empty for some $i,j\in\{1,\dots,n\}$, the formulas $\phi$ and $\phi'$ are equivalent.
	But $\phi'$ contains less non-Horn clauses than $\phi$, contradicting the choice of $\phi$.
	If there are distinct $i,j\in\{1,\dots,n\}$ and an integer $p\in E_{i,j}$ such that the formula $\phi\land\suc^p(x_i,x_j)$ is not equivalent
	to a Horn formula we reach a contradiction, as the relation defined by $\exists x_j(\phi\land \suc^p(x_i,x_j))$ is not Horn and has a smaller arity.
	Therefore, for each distinct $i,j$ and $p\in E_{i,j}$, the formula $\chi:=\phi\land\suc^p(x_i,x_j)$ is equivalent to a Horn formula,
	and by Lemma~\ref{lem:Horn-iff-reduced-Horn}, it even \emph{reduces} to a Horn formula.
	Note that since $\psi$ is a positive clause, the only way to make $\chi$ Horn by removing
	literals is to remove all but one literals in $\psi$.
	This means that there exists a literal $\ell^{i,j}_p$ of $\psi$ such that
	$$\phi\land \suc^p(x_i,x_j)\models \ell^{i,j}_p,$$
	and in this case we could reduce $\phi$ by replacing $\psi$ by $\ell^{i,j}_p$ and have an equivalent formula.
	Let $q$ be the qe-degree of $\phi$.
	If $p$ is in $E_{i,j}$ and is greater than $nq$, then we may take $\ell^{i,j}_p$
	to be $\ell^{i,j}_{nq+1}$, by the substitution lemma.
	
	Assume first that for some distinct $i,j$, $E_{i,j}$ is finite.
	Then $\phi$ is equivalent over $(\Z;\suc)$ to the formula
	$$\chi:=\phi' \land \bigwedge_{p\in E_{i,j}} (\suc^p(x_i,x_j) \Rightarrow \ell^{i,j}_p) $$
	which has fewer non-Horn clauses than $\phi$.
	Indeed, $\phi$ implies $\chi$ directly from the hypotheses we have. Conversely, if $s$ satisfies $\chi$
	one of two cases occur. Either some $\ell^{i,j}_p$ if satisfied by $s$,
	and then $s$ satisfies $\psi$ and $\phi$.
	Or we must have $s(x_j)\neq s(x_i)+p$ for every $p\in E_{i,j}$, i.e., $s(x_j)-s(x_i)\not\in E_{i,j}$.
	Since $s$ is known to satisfy $\phi'$, by definition of $E_{i,j}$ it must also satisfy $\psi$, whence we get
	that $s$ satisfies $\phi$. Note that $\chi$ contains fewer non-Horn clauses than $\phi$,
	which contradicts the choice of $\phi$.
	
	Therefore all the $E_{i,j}$ are cofinite for distinct $i$ and $j$, and therefore $nq+1\in E_{i,j}$ for every distinct $i,j$.
	As a consequence, if $s$ is a satisfying assignment for $\phi$ and $i,j\in\{1,\dots,n\}$ are distinct and such that
	$|s(x_i)-s(x_j)|>nq$, $s$ must satisfy the literal $\ell^{i,j}_{nq+1}$. 
	Let $s$ be an assignment of $\phi$ such that for every distinct $i,j\in\{1,\dots,n-1\}$, we have $|s(x_i)-s(x_j)|>2(n+1)q$  (such an assignment exists,
	by the fact that the $(n-1)$-projection of $R$ onto $\{1,\dots,n-1\}$ is Horn, and since the binary projections of $R$ all have cofinite distance degree).
	For $\psi$ to be satisfied by $s$, there must exist an $i\in\{1,\dots,n-1\}$ such that $|s(x_i)-s(x_n)|\leq q$
	(remember that $\psi$ only contains positive literals of degree at most $q$).
	Let $k\in\{1,\dots,n-1\}$ be different than $i$.
	Note that $|s(x_k)-s(x_i)|>nq$ and $|s(x_k)-s(x_n)|>nq$.
	Then the literal $\ell^{k,i}_{nq+1}$ relates $x_i$ and $x_n$,
	and so does the literal $\ell^{k,n}_{nq+1}$, because $x_i$ and $x_n$
	are the only variables that are able to satisfy a positive literal.
	Let now $t$ be an assignment of $\phi$ such that $|t(x_i)-t(x_n)|>2(n+1)q$.
	Either $|t(x_k)-t(x_i)| > nq$ or $|t(x_k)-t(x_n)| > nq$.
	In the first case, $\ell^{k,i}_{nq+1}$ must be satisfied by $t$.
	But $\ell^{k,i}_{nq+1}$ is a literal of the form $\suc^p(x_n,x_i)$ with $|p|\leq q$, and
	$|t(x_i)-t(x_n)|>nq$, so $t$ cannot satisfy $\ell^{k,i}_{nq+1}$. Similarly, in the second case,
	$t$ must satisfy $\ell^{k,n}_{nq+1}$, which is impossible
	since this literal is of the form $\suc^p(x_n,x_i)$ for $|p|\leq q$.
	We have reached a contradiction.
	Therefore, we must have $n=2$.
\end{proof}

Note that a binary relation that is not Horn is non-trivial and has finite distance degree.

\begin{lemma}\label{lem:not-positive-not-Horn-hard}
	Let $\Gamma$ be a non-positive expansion of $(\Z;\suc)$ such that
	$\Q.\Gamma$ omits $\si$ as a polymorphism.
	Then $\Csp(\Gamma)$ is NP-hard.
\end{lemma}
\begin{proof}
	By Proposition~\ref{prop:non-Horn-binary-finite}, there exists a non-trivial binary relation $T$ of finite distance degree with a pp-definition in $\Gamma$.
	By Lemma~\ref{lem:positive-iff-endo-collapses}, there exists a 
	non-trivial binary relation $N$ pp-definable in $\Gamma$ and which has infinite distance degree.
	The relation defined by $N(x,y)\land N(y,x)$ in $\Gamma$ is symmetric and has infinite distance degree,
	and is again pp-definable in $\Gamma$, so we will assume that $N$ is already symmetric.
	Let $a$ be the smallest positive integer such that $(0,b)$ is in $N$ for all $b\geq a$.
	With $\suc$ and pp-definition, we may assume that $T$ contains $(0,0)$
	and pairs $(0,b),(0,2b)$ with $b\geq a$.
	Let $G$ be the undirected graph whose vertices are the integers $v$ such that $(0,v)\in T$,
	and where $v$ and $w$ are adjacent if $(v,w)\in N$.
	This graph contains the triangle $(0,b),(b,2b),(0,2b)$, so that $G$
	is not bipartite and $\Csp(G)$ is NP-hard by~\cite{HellNesetril}.
	Furthermore, $\Csp(G)$ is polynomial-time reducible to $\Csp(\Gamma)$:
	if $\exists x_1,\dots,x_n \phi$ is an instance of $\Csp(G)$,
	create an instance of $\Csp(\Gamma)$ by adding an existentially quantified
	variable $z$, and by adding the constraints $T(z,x_i)$ for all $i$.
	It is clear that this instance is satisfiable iff the original instance is satisfiable
	in $G$, using the fact that the automorphism group of $\Gamma$ is transitive.
	This proves that $\Csp(\Gamma)$ is NP-hard.
\end{proof}

\subsection{The positive case}

We prove in this section that a positive first-order expansion $\Gamma$ of $(\Z;\suc)$ which is not preserved by any $d$-modular maximum or minimum
has an NP-hard CSP. As in the non-positive case and Proposition~\ref{prop:non-Horn-binary-finite}, an important step of the classification
is to show that there exists a non-trivial binary relation with a pp-definition in $\Gamma$. 

Let $R$ be a relation of arity $n$ with a first-order definition $\phi$ over a structure $\Gamma$. We say that $R$ is \emph{$r$-decomposable}
if $\phi(x_1,\dots,x_n)$ is equivalent in $\Gamma$ to $$\bigwedge_J \exists_{j\not\in J} x_j  \phi(x_1,\dots,x_n)$$
where $J$ ranges over all the $r$-element subsets of $\{1,\dots,n\}$.
The following lemma is a positive equivalent to Lemma~38 in~\cite{BodDalMarMotPin}, and its proof
is essentially the same.
Intuitively this is because in both cases the binary relations that
are pp-definable in $\Gamma$ have either a finite distance degree or are $\Z^2$
(if $\Gamma$ has finite distance degree this is immediate, and when $\Gamma$ is positively definable in $(\Z;\suc)$ this is the content of Lemma~\ref{lem:positive-iff-endo-collapses}).
For the sake of completeness, we reproduce the proof with the necessary adjustments.

\begin{lemma}\label{lem:positive-not-modular-binary-relation}
	Let $\Gamma$ be a positive fo-expansion of $(\Z;\suc)$ that does not admit a modular max or modular min  polymorphism.
	Then there is a non-trivial binary relation pp-definable in $\Gamma$ which has a finite distance degree.
\end{lemma}
\begin{proof}
	The binary relations pp-definable in $\Gamma$ are either trivial, or non-trivial and have a finite distance degree,
	by the fact that $\Gamma$ is positive and Lemma~\ref{lem:positive-iff-endo-collapses}.
	Suppose for contradiction that all the binary relations with a pp-definition in $\Gamma$ are trivial.
	
	If every relation $S$ pp-definable in $\Gamma$ were $2$-decomposable
	then $S$ would be invariant under a modular max or modular min operation:
	indeed, we assumed that the binary relations pp-definable in $\Gamma$
	are already pp-definable in $(\Z;\suc)$, so that the $2$-decomposable
	relations that have a pp-definition in $\Gamma$ already have a pp-definition
	in $(\Z;\suc)$, which means that they are
	preserved by the max$_d$ for all $d\geq 1$.
	Hence, there is a relation $S$ pp-definable in $\Gamma$ that is not $2$-decomposable.
	This implies that, by projecting out coordinates from $S$,
	we can obtain a relation $R$ of arity $r \geq 3$ which
	is not $(r-1)$-decomposable.

	This implies, in particular, that there exists a tuple $(a_1,\ldots,a_r) \notin R$
	such that for all $i \in \{1,\dots,r\}$, $(a_1,\ldots,p_i,\ldots,a_r) \in R$
	for some integer $p_i$.
	By replacing $R$ by the relation with the pp-definition
	\[ \exists y_1,\ldots, y_r \, \big(\bigwedge_{i \in [r]} (y_i = x_i + a_i) \wedge R(y_1,\ldots, y_r)\big) \]
	we can further assume that $a_i = 0$ for all $i \in [r]$.
	We can also assume, w.l.o.g., that $p_1\neq -p_2$ because $r\geq 3$.

	Suppose that the arity of $R$ is greater than $3$,
	and consider now the ternary relation $T(x_1,x_2,x_3)$
	defined by $R(x_1,x_2,x_3,\ldots,x_3)$.
	Suppose there is a $z$ so that $R(0,0,z,\ldots,z)$,
	then $T$ would not be $2$-decomposable since $(0,0,0)\not\in T$, although
	$(p_1,0,0),(0,p_2,0)$, and $(0,0,z)$ are all in $T$, which contradicts the
	minimality of the arity of $R$.
	If there is no such $z$ then
	$\exists x_3 R(x_1,x_2,x_3,\ldots,x_3)$
	defines a binary relation omitting $(0,0)$ and containing $(0,-p_1)$ and $(0,p_2)$.
	This relation is non-trivial (and hence has infinite distance degree, by the positivity assumption)
	and binary, contradiction.

	Thus we are in the situation in which $r=3$.
	Note that every binary projection of $R$ is $\Z^2$:
	otherwise one such binary projection (w.l.o.g\ $\exists x_1  R(x_1,x_2,x_3)$)
	would be of the form $x_3=x_2+p$ for some $p\in\Z$.
	Let $(a,b,c)$ be such that $(a,b)$ is in the projection of $R$ onto $\{1,2\}$,
	$(a,c)$ is in the projection of $R$ onto $\{1,3\}$,
	and $(b,c)$ is in the projection of $R$ onto $\{2,3\}$ (i.e. $c=b+p$).
	Since $(a,b)$ is in the first projection of $R$, there exists $d\in\Z$
	such that $(a,b,d)$ is in $R$, but since the third projection is trivial
	we have $d=b+p=c$, so that $(a,b,c)$ is in $R$ and $R$ is $2$-decomposable,
	contradicting our assumptions. Thus every binary projection of $R$ is $\Z^2$.
	
	Let $\phi(x_1,x_2,x_3)$ be a positive formula in reduced DNF,
	and let $R$ be the relation defined by $\phi$ over $(\Z;\suc)$.
	This formula has at least two disjuncts, otherwise $R$ would be pp-definable
	over $(\Z;\suc)$.
	Each disjunct contains at most two literals, because it suffices
	to describe only two distances between three variables to determine the type of a triple of integers.
	We claim that there is a disjunct in $\phi$ that consists of only one literal.
	If that was not the case, every disjunct $\mathcal D_i$ would have two literals
	and would be equivalent	to $\suc^{p_i}(x_2,x_1)\land \suc^{q_i}(x_3,x_1)$ for some $p_i,q_i\in\Z$.
	In this case, the formula $\exists x_2.  \phi(x_1,x_2,x_3)$ defines a binary relation
	with finite distance degree, contradicting the last sentence in the previous paragraph.
	Furthermore, there are at least two such disjuncts:
	if there is only one, say $\suc^p(x_2,x_1)$, the relation defined
	by $\exists x_3.  \phi(x_1,x_2,x_3)$ is binary and has a finite distance degree, a contradiction.
	Hence there are at least two disjuncts in $\phi$ that contain only one literal.
	One of $x_1,x_2,x_3$ must appear twice in those literals, 
	and we may assume by permuting the variables that it is $x_1$.
	Let us write these literals as $\suc^p(x_2,x_1)$ and $\suc^q(x_3,x_1)$, for $p,q\in\Z$.
	Then the formula $\exists x_3 \, \big (\phi(x_1,x_2,x_3)\land \suc^{p-q+1}(x_2,x_3)\big )$
	is equivalent to a binary DNF which is reduced and contains the two disjuncts
	$\suc^p(x_2,x_1)$ and $\suc^{p+1}(x_2,x_1)$.
	The relation defined by this formula has finite distance degree,
	again contradicting our assumptions.

	It follows that that there exists a non-trivial binary relation pp-definable in $\Gamma$,
	and this relation has finite distance degree by positivity of $\Gamma$.
\end{proof}

\begin{definition}\label{def:progression}
A \emph{$d$-progression} is a set of the form
$[a,b\mid d] := \{a,a+d,a+2d,\ldots,b\}$, for $a \leq b$ with $b - a$ divisible by $d$.
\end{definition}

\begin{lemma}[Lemma~43 in~\cite{BodDalMarMotPin}]\label{lem:generatedprogression}
 Let $S \subset \Z$ be finite with $|S|>1$,  and let $d$ be the greatest common divisor of all $a-a'$ for $a,a'\in S$. Then
for any $d$-progression $T$, the relation $\Diff{T}$ 
is pp-definable in $(\Z; \suc,  \Diff{S})$.   
\end{lemma}

\begin{proposition}
\label{prop:new-with-1-progression}
Let $\Gamma$ be a first-order 
expansion of $(\Z;\suc)$, and 
$S \subset \Z$ a $1$-progression, $|S| > 1$, such that $\Diff{S}$ is pp-definable 
in $\Gamma$. 
Then $\Gamma$ is preserved by max or min;
or $\Csp(\Gamma)$ is NP-hard. 
\end{proposition}
\begin{proof}
    Suppose that $\Gamma$ is not preserved by $\max$ nor $\min$.
    Therefore, there exist in $\Gamma$ a relation $R\subseteq\Z^n$ that is not preserved by $\max$ and a relation $T\subseteq\Z^m$ which is not preserved by $\min$.
    This means that there are tuples $\tuple a,\tuple b$ in $R$ such that $\max(\tuple a,\tuple b)$ is not in $R$
    and similarly for $T$.
    By hypothesis and Lemma~\ref{lem:generatedprogression}, all the 1-progressions are definable in $\Gamma$.
    Let $M$ be $\max_{i,j}\{|a_i-a_j|,|b_i-b_j|\}$, so that $\Dist{[0,M]}$ is pp-definable in $\Gamma$.
    Define the relation $R^*$ by $\exists x_1,\dots,x_n(R(x_1,\dots,x_n) \land \bigwedge_{i\neq j}(x_i,x_j)\in\Dist{[0,M]})$
    and similarly define $T^*$. We have that $\tuple a,\tuple b$ are in $R^*$ by construction,
    and still $\max(\tuple a,\tuple b)\not\in R^*$ since $R^*\subseteq R$.
    Moreover, $T^*$ and $R^*$ have finite distance degree.
    By Proposition~47 in~\cite{BodDalMarMotPin}, $\Csp(\Z;\suc,\Diff{S},R^*,T^*)$ is NP-hard, 
    therefore $\Csp(\Gamma)$ is also NP-hard.
\end{proof}

\subsection{Concluding the Classification}

\begin{thm:succ}
Let $\Gamma$ be a first-order expansion of $(\Z;\suc)$. 
Then at least one of the following is true:
\begin{enumerate}
	\item $\Gamma$ is positive and preserved by a modular max or a modular min,
	\item $\Gamma$ is non-positive and $\Q.\Gamma$ is preserved by $\si$,
	\item $\Csp(\Gamma)$ is NP-hard.
\end{enumerate}
\end{thm:succ}
\begin{proof}
	Suppose first that $\Gamma$ is non-positive.
	If $\Q.\Gamma$ is preserved by $\si$, we are done.
	Otherwise, 
	Lemma~\ref{lem:not-positive-not-Horn-hard}
	states that $\Csp(\Gamma)$ is NP-hard.
	If $\Q.\Gamma$ is positive, suppose that $\Gamma$
	 omits all modular max and min polymorphisms. 
	By Lemma~\ref{lem:positive-not-modular-binary-relation} 
	there exists a non-trivial binary relation $R$ with a finite distance degree with a pp definition in $\Gamma$.
If $R$ is not a $d$-progression for any $d\geq 1$,
then $\Csp(\Gamma)$ is NP-hard by Lemma~44 in~\cite{BodDalMarMotPin}.

Finally, if $R$ is a non-trivial $d$-progression and $\Gamma$ is not preserved by $\max_d$ or $\min_d$, then $\Gamma/d$ is not preserved by $\max$ or $\min$.
Moreover, $\Gamma$ pp-defines a non-trivial $d$-progression so $\Gamma/d$ pp-defines a non-trivial $1$-progression, which means that Proposition~\ref{prop:new-with-1-progression} applies
and that CSP$(\Gamma/d)$ is NP-hard.
Now we reduce  CSP$(\Gamma/d)$ to CSP$(\Gamma)$ to prove that the latter is also NP-hard.
Let $q$ be the qe-degree of $\Gamma$ and note that an instance of $\Gamma$ on $n$ variables has a solution iff it has a solution on the interval $[0,qn]$.
From an instance $\Phi$ of CSP$(\Gamma/d)$ we build an instance $\Psi$ of CSP$(\Gamma)$.
To build $\Psi$ from $\Phi$, we augment with a new variable $z$ as well as $qn$ new variables $x_1\ldots,x_{qn}$ for each extant variable $x$ of $\Psi$.
Then $\Psi$ is as $\Phi$ but with the additional constraints $\mathrm{Dist}_{[0,qdn|d]}(x,z)$,
where we define $\mathrm{Dist}_{[0,qd(n+1)|d]}(x,z)$ by $\mathrm{Dist}_{[0,d|d]}(x,x_1) \wedge \mathrm{Dist}_{[0,d|d]}(x_1,x_2) \wedge \ldots \wedge \mathrm{Dist}_{[0,d|d]}(x_{qn},z)$.
It is straightforward to see that  $\Gamma/d \models \Phi$ iff $\Gamma \models \Psi$ and the result follows.  
\end{proof}

\begin{thm:main}
Let $\Gamma$ be a reduct of $(\Z;<)$ with finite signature.
Then there exists a structure $\Delta$ such that $\Csp(\Delta)$ equals
$\Csp(\Gamma)$ and one of the following cases applies. 
\begin{enumerate}
\item $\Delta$ has a finite domain, and the CSP for $\Gamma$ is conjectured to be in Ptime or NP-complete~\cite{FederVardi}.
\item $\Delta$ is a reduct of $(\Q;<)$, and the complexity of $\Csp(\Delta)$ has been classified in~\cite{tcsps-journal}. 
\item $\Delta$ is a reduct of $(\Z;<)$ and preserved by a modular max or modular min.
In this case, $\Csp(\Gamma)$ is in Ptime.
\item $\Delta$ is a reduct of $(\Z;\suc)$ that is preserved by a binary injective function preserving $\suc$. 
In this case, $\Csp(\Gamma)$ is in Ptime. 
\item $\Csp(\Gamma)$ is NP-complete. 
\end{enumerate}
\end{thm:main}
\begin{proof}
Let $\Gamma$ be a finite signature reduct
of $(\Z;<)$. By Proposition~\ref{prop:inNP}, $\Csp(\Gamma)$ is in NP. If $\Gamma$ is homomorphically equivalent to a finite structure, we are in case one of the statement and there is nothing to be shown. Otherwise, 
Theorem~\ref{thm:class} implies that there
exists a reduct $\Delta$ of $(\Z;<)$ such that $\Csp(\Gamma)$ equals $\Csp(\Delta)$,
and one of the following cases applies.
\begin{enumerate}
\item $\Delta$ is a reduct of $(\Q;<)$;
by Theorem 50 in \cite{tcsps-journal}, $\Csp(\Delta)$
is in P or NP-complete, we are in case 2 of the statement.
\item For all $k \geq 1$, the relation 
$\Dist{\{k\}}$ 
is pp-definable; in this case, $\Csp(\Gamma)$ and $\Csp(\Delta)$ are NP-hard by Proposition~\ref{prop:26}. Hence, we are in case four of the statement. 
\item The relation $\suc$ is pp-definable in $\Delta$.
Suppose that neither item 3 nor item 4 applies.
If $<$ is pp-definable in $\Delta$, then 
Lemma~\ref{lem:not-max-NPhard} implies that $\Csp(\Delta)$ is NP-hard,
and we are in case five of the statement.
Otherwise $\Delta$ is a reduct of $(\Z;\suc)$, by Lemma~\ref{lem:define-order}.
In this case, Theorem~\ref{thm:succ} implies that one of the following cases applies. 
\begin{itemize} 
\item $\Delta$ is a positive expansion of $(\Z;\suc)$ and is preserved by a modular max or modular min polymorphism.
In this case, $\Csp(\Gamma)$ is in Ptime by Theorem~\ref{thm:tractability-modular-semilattice}, and we are in case three of the statement. 
\item $\Q.\Delta$ is a non-positive expansion of $(\Z;\suc)$ preserved by $\si$. 
In this case, $\Csp(\Gamma)$ is in Ptime by Corollary~\ref{cor:si-tract}, and we are in case four of the statement.
\item $\Csp(\Delta)$ is NP-hard. We are in case five of the statement. 
\end{itemize}
\end{enumerate}
\end{proof}

\section{Open Problems and Future Work}
\label{sect:discussion}
In this article, the complexity of 
$\Csp(\Gamma)$ has been classified for all relational structures $\Gamma$ over the integers where the relations are first-order definable over $({\mathbb Z};<)$, assuming
the Feder-Vardi conjecture. 
This class of CSPs subsumes the class of \emph{temporal CSPs}~\cite{tcsps-journal} and the class of \emph{distance CSPs} where the constraints are first-order definable over the integers with the successor relation~\cite{BodDalMarMotPin}. 

These results are important foundations for the future investigation of the important class of CSPs where the constraint relations are
definable in Presburger arithmetic, i.e., definable over $(\Z;+,<)$. We give here two possible classification projects that can improve our understanding 
of the complexity of problems expressible in Presburger arithmetic.
\begin{itemize}
\item Note that every integer has a first-order definition in $(\Z;+,<)$. The same is true for the comparatively simpler structure $(\Z;\suc,0)$.
While our results imply a complexity dichotomy for the CSPs of reducts of $(\Z;\suc)$ (at least, those that are not the CSP of a finite structure),
the techniques we employed cannot handle the case of $(\Z;\suc,0)$. The principal reason is that we used extensively the transitivity of the automorphism group of $(\Z;<)$, while $(\Z;\suc,0)$ is \emph{rigid}, i.e., has no automorphisms beside the identity function. 

On the other hand,
classifying the complexity of CSPs that are first-order definable with infinitely many constants can be reduced to proving the algebraic dichotomy conjecture for finite-domain CSPs~\cite{LiftingTractability}. 
It is therefore an interesting question whether the two results can be combined to obtain a complexity classification for the reducts of $(\Z;\suc,0)$. 
\item 
Secondly, a result for Presburger arithmetic would in particular give a complexity classification for the CSPs of reducts of $(\Z;+)$ or reducts of $(\Z;+,1)$.
The latter endeavour has been started by the authors and Marcello Mamino~\cite{ExpansionsPlusConstants},
with a complexity classification of the CSPs of the reducts of $(\Z;+,1)$ that contain $+$ in their signature.
\end{itemize}

An important family of open problems concerns the complexity of constraint satisfaction problems over the integers 
for \emph{infinite relational signatures}. 
When working with infinite signatures, we have to specify how the relation symbols are represented in the input instances. If we represent the relation symbol for a relation $R$ by a quantifier-free definition of $R$ using
atomic formulas of the form $y \leq x+ c$ where $c \in {\mathbb Z}$ is represented in binary, one can 
formulate the famous Max-Atoms problems as a CSP in this class. 
The Max-Atoms problem is known to be polynomial-time equivalent to determining the winner in Mean Payoff Games~\cite{and-or-scheduling,Max-atoms}, which is a problem in the intersection of NP and coNP, but not known to be in P. Our proofs make crucial use of finite signatures and bounded quantifier-degrees; but 
many of the statements in this article could hold for reducts of $({\mathbb Z};<)$ 
with \emph{infinite} relational signature
(we could not find counterexamples). 

\bibliographystyle{acm}
\bibliography{local}

\end{document}